\documentclass[a4paper,12pt,oneside]{amsart}

\usepackage{amsmath,amsfonts,amssymb,amsthm,amscd,latexsym,euscript,xypic}
\usepackage[all]{xy}

\usepackage{mathtools}

\usepackage{pdfsync}

\newtheorem{theorem}{Theorem}[section]

\newtheorem{claim}[theorem]{Claim}
\newtheorem{corollary}[theorem]{Corollary}

\newtheorem{definition}[theorem]{Definition}

\newtheorem{proposition}[theorem]{Proposition}
\newtheorem{prop}[theorem]{Proposition}

\newcommand{\Ae}{{\mathcal A}}

\newcommand{\Ce}{{\mathcal C}}
\newcommand{\De}{{\mathcal D}}
\newcommand{\Ee}{{\mathcal E}}

\newcommand{\He}{{\mathcal H}}

\renewcommand{\emptyset}{\varnothing}

\numberwithin{equation}{section}

\newcommand{\WVP}{\rm{WVP}}
\newcommand{\VP}{\rm{VP}}
\newcommand{\SWVP}{\rm{SWVP}}
\newcommand{\PRP}{\rm{PRP}}

\newcommand{\crit}{{\rm crit}}

\newcommand{\OR}{\mathrm{OR}}

\begin{document}

\title[The $\WVP$ for definable classes of structures]{The Weak Vop\v enka Principle for definable classes of structures}
\author{Joan Bagaria}
\address{ICREA (Instituci\'o Catalana de Recerca i Estudis Avan\c{c}ats) and
\newline \indent Departament de Matem\`atiques i Inform\`atica, Universitat de Barcelona. 
Gran Via de les Corts Catalanes, 585,
08007 Barcelona, Catalonia.}
\email{joan.bagaria@icrea.cat}

\author{Trevor M. Wilson}
\address{ Department of Mathematics
Miami University,
236 Bachelor Hall, 
Oxford, OH 45056}

\email{twilson@miamioh.edu }
\thanks{Part of this research was  supported by  the Generalitat de Catalunya (Catalan Government) under grant SGR 270-2017, and by the Spanish
Government under grant MTM2017-86777-P}

\date{\today }

\begin{abstract}
We give a level-by-level analysis of the Weak  Vop\v{e}nka Principle  for definable classes of relational structures ($\WVP$), in accordance with the complexity of their definition, and we determine the large-cardinal strength of each level. Thus, in particular we show that $\WVP$ for $\Sigma_2$-definable classes is equivalent to the existence of a strong cardinal. The main theorem (\ref{thm2}) shows, more generally, that $\WVP$ for $\Sigma_n$-definable classes  is equivalent to the existence of a $\Sigma_n$-strong cardinal (definition \ref{defSigmaStrong}). Hence, $\WVP$ is equivalent to the existence of a $\Sigma_n$-strong cardinal, all $n<\omega$.
\end{abstract}
\maketitle

\section{Introduction}

The Vop\v{e}nka Principle ($\VP$), which asserts that there is no rigid proper class of graphs, is a well-known strong large-cardinal principle\footnote{As related in \cite{AdR:LPAC}, \emph{``The story of Vop\v{e}nka's principle [...] is that of a practical joke which misfired: In the 1960's P. Vop\v{e}nka was repelled by the multitude of large cardinals which emerged in set theory. When he constructed, in collaboration with Z. Hedrl\'{\i}n and A. Pultr, a rigid graph on every set [...], he came to the conclusion that, with some more effort, a large rigid class of graphs must surely be also constructible. He then decided to tease set-theorists: he introduced a new principle (known today as Vop\v{e}nka's principle), and proved some consequences concerning large cardinals. He hoped that some set-theorists would continue this line of research (which they did) until somebody showed that the principle was nonsense. However the latter never materialized -- after a number of unsuccessful attempts at constructing a large rigid class of graphs, Vop\v{e}nka's principle received its name from Vop\v{e}nka's disciples."}}  (see \cite{Kan:THI}). Properly formulated as a first-order assertion, $\VP$ is a schema, that is, an infinite collection of statements, one for each first-order formula defining a proper class of graphs and asserting that the class defined by the formula is not rigid, i.e., there is some non-identity morphism. An equivalent formulation of $\VP$ as a first-order schema is given by restricting $\VP$ to proper classes of graphs that are definable with a certain degree of complexity, according to the Levy hierarchy of formulas $\Sigma_n$, $n<\omega$ (see \cite{Jech}). Thus, writing $\mathbf{\Sigma_n}$-$\VP$ for the first-order assertion that every $\Sigma_n$-definable (with parameters), proper class of graphs is rigid, we have that $\VP$ is equivalent to $\mathbf{\Sigma_n}$-$\VP$, for every $n$. As shown in \cite{Ba:CC, BCMR}, $\mathbf{\Sigma_1}$-$\VP$ is provable in ZFC, while $\mathbf{\Sigma_2}$-$\VP$ is equivalent to the existence of a proper class of supercompact cardinals, $\mathbf{\Sigma_3}$-$\VP$ is equivalent to the existence of a proper class of extendible cardinals, and $\mathbf{\Sigma_n}$-$\VP$ is equivalent to the existence of a proper class of $C^{(n-2)}$-extendible cardinals,  for $n\geq 3$. The level-by-level analysis of  $\mathbf{\Sigma_n}$-$\VP$, $n<\omega$, and the determination of their large-cardinal strength yielded new results in category theory, homology theory, homotopy theory, and universal algebra (see \cite{BCMR}). 
For example,  the existence of cohomological localizations in the homotopy category of simplicial sets (Bousfield Conjecture) follows from  $\mathbf{\Sigma_2}$-$\VP$.

The role of $\VP$ in category theory  has a rich history. The first equivalences of $\VP$ with various category-theoretic statements were announced by  E. R. Fisher  in \cite{F:VP}. Further equivalences were proved over the next two decades by Ad\'amek, Rosick\'y, Trnkov\'a, and others. Their work showed that  under $\VP$ \emph{``the structure of locally presentable categories becomes much more transparent"} (\cite{AdR:LPAC} p. 241). For example, the statement that a category is locally presentable if and only if it is complete and bounded is equivalent to $\VP$. And so is the statement that every orthogonality class in a locally presentable category is a small-orthogonality class (\cite{RTA:UPLPC},  \cite{AdR:LPAC} 6.9, 6.14). Of the many category-theoretic statements now known to be equivalent to $\VP$, the following one (see \cite{AdR:LPAC} 6.D) turned out to be of particular interest:
\begin{itemize}
\item[(1)] Every full subcategory of a locally presentable category $\mathcal{K}$ closed under colimits is coreflective in $\mathcal{K}$.
\end{itemize}
What made  (1) particularly interesting is that its dual statement 
\begin{itemize}
\item[(2)] Every full subcategory of a locally presentable category $\mathcal{K}$ closed under limits is reflective in $\mathcal{K}$.
\end{itemize}
while being a consequence of (1), could not be proved equivalent to it. Since   $\VP$ -- hence also (1) -- was known to be equivalent to
$$\mathbf{Ord}\mbox{ cannot be fully embedded into }\mathbf{Gra}$$
(see \cite{AdR:LPAC} 6.3), while statement (2) was proved equivalent to
$$\mathbf{Ord}^{op}\mbox{ cannot be fully embedded into }\mathbf{Gra}$$
(see \cite{AdR:LPAC} 6.22, 6.23), the latter assertion was then called \emph{The Weak Vop\v{e}nka Principle} (\WVP). The term \emph{Weak} was aptly given, for it is readily shown that $\VP$ implies $\WVP$ (\cite{ART}; Proposition \ref{propVPimpliesWVP} below). The question then remained if $\WVP$ implied $\VP$. Using a result of J. R. Isbell \cite{I:AS}, which showed that $\mathbf{Ord}^{op}$ is bounded iff there is no proper class of measurable cardinals, Ad\'amek-Rosick\'y \cite{AdR:LPAC}  proved that $\WVP$ implies the existence of a proper class of measurable cardinals. This was seen as a first step in showing that $\WVP$ was indeed a strong large-cardinal principle, perhaps even  equivalent to $\VP$. Much work was devoted to trying to obtain stronger large cardinals from it, e.g., strongly compact or supercompact cardinals, but to no avail. A further  natural   principle,  between $\VP$ and $\WVP$, called the \emph{Semi-Weak Vop\v{e}nka Principle} ($\SWVP$), was introduced in \cite{AdR:ILPC}  and the further question of the equivalence between the three principles, $\WVP$, $\SWVP$ and $\VP$, remained open. The problem was finally solved in 2019 by the second author of the present paper. In \cite{W:LCSWVP} he showed that $\WVP$ and $\SWVP$ are equivalent, and they are also equivalent to the large-cardinal principle ``$\OR$ is Woodin", whose consistency strength is known to be  well below the existence of a supercompact cardinal, thereby showing that $\WVP$ cannot imply $\VP$ (if consistent with ZFC). 

In the present paper we carry out a level-by-level analysis of $\WVP$ and $\SWVP$ similar to the analysis of $\VP$ done in \cite{Ba:CC, BCMR}. Thus, we prove  the equivalence of both $\mathbf{\Sigma_n}$-$\WVP$ and $\mathbf{\Sigma_n}$-$\SWVP$ (see definition \ref{defWVP} below) with large-cardinals, for every $n\geq 2$.   In particular, we show that $\mathbf{\Sigma_2}$-$\WVP$ and $\mathbf{\Sigma_2}$-$\SWVP$ are equivalent to the existence of a proper class of strong cardinals. The main theorems (\ref{thm2}, \ref{thm3}) show, more generally, that $\mathbf{\Sigma_n}$-$\WVP$ and $\mathbf{\Sigma_n}$-$\SWVP$ are equivalent to the existence of a proper class of $\Sigma_n$-strong cardinal (definition \ref{defSigmaStrong}). Moreover, $\WVP$ and $\SWVP$ are equivalent to the existence of a $\Sigma_n$-strong cardinal, for every $n<\omega$. Our arguments   yield also a new proof of the second author's result from  \cite{W:LCSWVP} that $\WVP$ implies ``$\OR$ is Woodin" (corollary \ref{maincoro} below). The main difference between the two proofs is that while in the present paper we derive the  extenders witnessing ``$\OR$ is Woodin" from homomorphisms on products of relational structures with universe of the form $V_\alpha$, the proof in \cite{W:LCSWVP} uses homomorphisms of so-called $\mathcal{P}$-structures. We think, however, that it should be possible to do a similar level-by-level analysis as done here by using $\mathcal{P}$-structures instead.  A number of consequences in category theory should follow from our results. For instance, the statement that every $\mathbf{\Sigma_2}$-definable full subcategory of a locally-presentable category $\mathcal{K}$ closed in $\mathcal{K}$ under limits is reflective in $\mathcal{K}$, should be equivalent to the existence of a proper class of strong cardinals. See  \cite{AdR:LPAC} Chapter 6 for more examples.

\section{Preliminaries}

Recall that a \emph{graph} is a structure $G=\langle G, E_G\rangle$, where $G$ is a non-empty set and $E_G$ is a binary relation on $G$. If $G=\langle G, E_G\rangle$ and $H=\langle H, E_H\rangle$ are graphs, a map $h:G\to H$ is a \emph{homomorphism} if it preserves the binary relation, meaning that for every $x,y\in G$, if $xE_G y$, then $h(x)E_H h(y)$. 

A class $\mathcal{G}$ of graphs is called \emph{rigid} if there are no non-trivial homomorphisms between graphs in $\mathcal{G}$, i.e., the only homomorphisms are the identity morphisms $G\to G$, for $G\in \mathcal{G}$.

The original formulation of \emph{Vop\v{e}nka Principle} ($\VP$) (P. Vop\v{e}nka, ca. ~1960) asserts that there is no rigid proper class of graphs.  
As shown in \cite{AdR:LPAC} 6.A,  $\VP$
 is equivalent to the statement that the category $\mathbf{Ord}$ of ordinals cannot be fully embedded into the category $\mathbf{Gra}$ of graphs. That is, there is no sequence $\langle G_\alpha :\alpha \in \OR\rangle$ of graphs  such that for every $\alpha \leq \beta$ there exists exactly one homomorphism $G_\alpha \to G_\beta$, and no homomorphism $G_\beta \to G_\alpha$ whenever $\alpha < \beta$.
 
 The \emph{Weak Vop\v{e}nka Principle} ($\WVP$) (first introduced in \cite{ART}) is the statement \emph{dual} to $\VP$, namely that the opposite category  of ordinals, $\mathbf{Ord}^{op}$, cannot be fully embedded into $\mathbf{Gra}$. That is, there is no sequence $\langle G_\alpha :\alpha \in \OR\rangle$ of graphs  such that for every $\alpha \leq \beta$ there exists exactly one homomorphism $G_\beta \to G_\alpha$, and no homomorphism $G_\alpha \to G_\beta$ whenever $\alpha < \beta$.
 
 The \emph{Semi-Weak Vop\v{e}nka Principle} ($\SWVP$) (\cite{AdR:ILPC}) asserts that  there is no sequence $\langle G_\alpha :\alpha \in \OR\rangle$ of graphs  such that for every $\alpha \leq \beta$ there exists some (not necessarily unique) homomorphism $G_\beta \to G_\alpha$, and no homomorphism $G_\alpha \to G_\beta$ whenever $\alpha < \beta$.

Clearly, $\SWVP$ implies $\WVP$. The second author showed in \cite{W:WVP} that $\SWVP$ is in fact equivalent to $\WVP$. As shown in \cite{ART},  $\WVP$ is a consequence of $\VP$, and the same argument also shows that $\VP$ implies $\SWVP$. In fact, the argument shows the following:

\begin{prop}
\label{propVPimpliesWVP}
$\VP$ implies  that for every sequence $\langle G_\alpha : \alpha \in \OR \rangle$ of graphs there exist $\alpha <\beta$  with a homomorphism $G_\alpha \to G_\beta$.
\end{prop}

\begin{proof}
Suppose $\langle G_\alpha : \alpha \in \OR \rangle$ is a sequence of graphs. Without loss of generality, if $\alpha < \beta$, then $G_\alpha$ and $G_\beta$ are not isomorphic. Since there are only set-many (as opposed to proper-class-many) non-isomorphic graphs of any given cardinality, there exists a proper class $C\subseteq \OR$ such that $|G_\alpha | <|G_\beta|$ whenever $\alpha <\beta$ are in $C$.  For each $\alpha\in C$, add to $G_\alpha =\langle G_\alpha , E_\alpha \rangle$ a rigid binary relation $S_\alpha$ (\cite{VPH}), as well as the non-identity relation $\not =$, and consider the structure $A_\alpha =\langle G_\alpha , E_\alpha , S_\alpha , \not =\rangle$, which can be easily seen as a graph. Since the cardinalities are strictly increasing, by the $\not =$ relation there cannot be any homomorphism $A_\beta \to A_\alpha$ with  $\alpha <\beta$.  Also, because of the rigid relation $S_\alpha$, the identity is the only homomorphism $A_\alpha \to A_\alpha$. Since by $\VP$ the class $\{ A_\alpha :\alpha \in C\}$ is not rigid, there must exist $\alpha <\beta$ with a homomorphism $A_\alpha \to A_\beta$, hence a homomorphism $G_\alpha \to G_\beta$. 
\end{proof}

The definitions of $\VP$, $WVP$ and $\SWVP$ given above quantify over arbitrary classes,  so they are not first-order. Thus, a proper study of these principles must be carried out in some adequate class theory, such as NBG. In particular, the proof of last proposition can only be formally given in such class theory. We shall  however be interested in the forthcoming in the first-order versions of $\VP$, $WVP$ and $\SWVP$, which require to restrict to definable classes. 

\subsection{The $\VP$, $WVP$ and $\SWVP$ for definable classes}
Each of $\VP$, $WVP$ and $\SWVP$ can be formulated in the first-order language of set theory as a definition schema, namely as an infinite list of definitions, one for every natural number $n$, as follows: 
\begin{definition}
\label{defWVP}
Let $n$ be a natural number, and let $P$ be a set or a proper class. 

The \emph{$\Sigma_n(P)$-Vop\v{e}nka Principle ($\Sigma_n(P)\mbox{-}\VP$ for short)} asserts that there is no $\Sigma_n$-definable, with parameters in $P$, sequence $\langle G_\alpha :\alpha \in \OR\rangle$ of graphs  such that for every $\alpha \leq \beta$ there exists exactly one homomorphism $G_\alpha \to G_\beta$, and no homomorphism $G_\beta \to G_\alpha$ whenever $\alpha <\beta$.

The \emph{$\Sigma_n(P)$-Weak Vop\v{e}nka Principle ($\Sigma_n(P)\mbox{-}\WVP$ for short)} asserts that there is no $\Sigma_n$-definable, with parameters in $P$, sequence $\langle G_\alpha :\alpha \in \OR\rangle$ of graphs  such that for every $\alpha \leq \beta$ there exists exactly one homomorphism $G_\beta \to G_\alpha$, and no homomorphism $G_\alpha \to G_\beta$ whenever $\alpha <\beta$.

The boldface versions \emph{${\mathbf \Sigma_n}$-$\VP$ and ${\mathbf \Sigma_n}\mbox{-}\WVP$} are defined as $\Sigma_n(V)\mbox{-}\VP$ and $\Sigma_n(V)\mbox{-}\WVP$ respectively, i.e., any set is allowed as a parameter in the definitions.

$\Pi_n(P)$-$\VP$ and $\Pi_n(P)\mbox{-}\WVP$,  as well as $\mathbf{\Pi_n}$-$\VP$ and $\mathbf{\Pi_n}$-$\WVP$, and the lightface (i.e., without parameters) versions $\Sigma_n{\rm{-}}\VP$, $\Sigma_n{\rm{-}}\WVP$ and  $\Pi_n{\rm{-}}\VP$, $\Pi_n{\rm{-}}\WVP$,  are defined similarly. 

The \emph{Vop\v{e}nka Principle ($\VP$)}  is the schema asserting that the $\mathbf{\Sigma_n}{\rm{-}}\VP$ holds for every $n$. And the \emph{Weak Vop\v{e}nka Principle ($\WVP$)}  is the schema asserting that the $\mathbf{\Sigma_n}{\rm{-}}\WVP$ holds  for every $n$.

\medskip

If instead of requiring that for $\alpha \leq \beta$ there is exactly one homomorphism $G_\beta \to G_\alpha$ we only require that there is at least one, then we obtain the  \emph{Semi-Weak Vop\v enka Principle ($\SWVP$)}, formulated as the first-order schema  $\mathbf{\Sigma_n}{\rm{-}}\SWVP$, $n<\omega$.
\end{definition}

It is well-known that the category of structures in any fixed (many-sorted, infinitary) relational language can be fully embedded into $\mathbf{Gra}$  (see \cite{AdR:LPAC} 2.65). Thus,  if in the original definitions of $\VP$, $\WVP$ and $\SWVP$ one replaces  ``graphs" by ``structures in a fixed (many-sorted, infinitary) relational language", one obtains equivalent notions.  The same is true for the first-order formulations of these principles, but some extra care is needed to ensure there is no increase in the complexity of the definitions. In particular, in the case of infinite language signatures, an extra parameter for the language signature $\tau$, as well as a parameter for a rigid binary relation on a binary signature associated to $\tau$, may be needed in the definition. Namely, suppose  $\Gamma$ is one of the definability classes $\Sigma_n$, $\Pi_n$,  with $n\geq 1$, $P$ is a set or a proper class, and $\Ce$ is a $\Gamma$-definable, with parameters in $P$, class of  (possibly many-sorted) relational structures in a language type $\tau$, i.e., $\tau=\langle R_\alpha :\alpha <\lambda\rangle$, where each $R_\alpha$ is an $n_\alpha$-ary relation symbol,  $n_\alpha$ being some ordinal, possibly infinite. As in  \cite{AdR:LPAC} 2.65, there is a $\Delta_1$-definable (i.e., both $\Sigma_1$-definable and $\Pi_1$-definable), using $\tau$ as a parameter, one-sorted  binary type $\tau'$ (meaning that all the relations are binary), and also a $\Gamma$-definable, with parameters in $P$ plus $\tau$ as an additional  parameter, full embedding of $\Ce$ into  the category $\mathbf{Rel} \, \tau'$ of $\tau'$-structures and homomorphisms. Furthermore, there is a $\Delta_1$-definable, using $\tau$ and a rigid binary relation $r$ on $\tau'$ as   parameters, full embedding of $\mathbf{Rel} \, \tau'$ into $\mathbf{Gra}$. Hence, there is a $\Gamma$-definable (with parameters in $P$, plus $\tau$ and $r$ as additional parameters) full embedding of $\Ce$ into $\mathbf{Gra}$. Therefore, in the definitions of $\Gamma$-$\VP$, $\Gamma$-$\WVP$, and $\Gamma$-$\SWVP$ we may replace ``graphs" by ``structures in a fixed (many-sorted, infinitary) relational language" and obtain equivalent principles, provided we allow for the additional parameters ($\tau$ and $r$) involved. Let us, however, stress the fact that in the case of finite $\tau$, or even if $\tau$ is countable infinite and definable without parameters (e.g., recursive), then no additional parameters are involved, and therefore the versions of  $\Gamma$-$\VP$, $\Gamma$-$\WVP$, and $\Gamma$-$\SWVP$  for graphs and for relational structures are equivalent.

\subsection{Strong cardinals} 
Recall  that a cardinal $\kappa$ is \emph{$\lambda$-strong}, where $\lambda$ is a cardinal greater than $\kappa$,  if there exists an elementary embedding $j: V\to M$, with $M$ transitive, with critical point $\kappa$, and with $V_\lambda$ contained in $M$. A cardinal $\kappa$ is \emph{strong} if it is $\lambda$-strong for every cardinal $\lambda >\kappa$.

If $\kappa$ is a strong cardinal, then for every cardinal $\lambda >\kappa$ there exists an elementary embedding $j: V\to M$, with $M$ transitive, critical point $\kappa$, $V_\lambda$ contained in $M$, and $j(\kappa)>\lambda$. Moreover, if $\kappa$ is strong, then $V_\kappa \preceq_{\Sigma_2}V$. (See \cite{Kan:THI}.)

It is well-known that the notion of strong cardinal can be formulated in terms of  extenders (see \cite{Kan:THI}, section  26). Namely, 

\begin{definition}
\label{defext}
Given a cardinal $\kappa$, and $\beta >\kappa$, a \emph{$(\kappa,\beta)$-extender} is a collection $\Ee:=\{ E_a:a\in [\beta]^{<\omega}\}$ such that 
\begin{enumerate}
\item Each $E_a$ is a $\kappa$-complete ultrafilter over $[\kappa ]^{|a|}$, and $E_{a}$ is not $\kappa^+$-complete for some $a$. 
\item For each $\xi <\kappa$, there is  some  $a$ with  $\{ s\in [\kappa]^{|a|}:\xi \in s\}\in E_a$.
\item \emph{Coherence:} If $a\subseteq b$ are in $[\beta]^{<\omega}$, with $b=\{ \alpha_1,\ldots, \alpha_n\}$ and $a=\{ \alpha_{i_1}, \ldots , \alpha_{i_n}\}$, and  $\pi_{ba}:[\kappa]^{|b|}\to [\kappa]^{|a|}$ is the map given by $\pi_{ba}(\{ \xi_1,\ldots ,\xi_n\})=\{ \xi_{i_1},\ldots ,\xi_{i_n}\}$, then 
$$X\in E_a \quad \mbox{ if and only if } \quad \{ s\in [\kappa]^{|b|}: \pi_{ba}(s)\in X\}\in E_b\, .$$
\item  \emph{Normality:} Whenever $a\in [\beta]^{<\omega}$ and $f:[\kappa]^{|a|} \to V$ are such that $\{ s\in [\kappa]^{|a|}: f(s)\in {\rm{max}}(s)\}\in E_a$, there is $b\in [\beta]^{<\omega}$   with $a\subseteq b$ such that
$$\{ s\in [\kappa]^{|b|}: f(\pi_{ba}(s))\in s\}\in E_{b}\, .$$
\item \emph{Well-foundedness:} Whenever $a_m\in [\beta]^{<\omega}$ and $X_m\in E_{a_m}$ for $m\in \omega$, there is a function $d:\bigcup_m a_m \to \kappa$ such that $d`` a_m \in X_m$ for every $m$.
\end{enumerate}

\end{definition}

\begin{prop}
\label{prop-1}
A cardinal $\kappa$ is $\lambda$-strong if and only if there exists a  $(\kappa, |V_{\lambda}|^+)$-extender $\mathcal{E}$ such that $V_{\lambda}\subseteq \overline{M}_{\mathcal{E}}$ and $\lambda <j_{\mathcal{E}}(\kappa)$ (where $\overline{M}_{\mathcal{E}}$ is the transitive collapse of the direct limit ultrapower $M_\mathcal{E}$ of $V$ by $\mathcal{E}$, and $j_{\mathcal{E}}:V\to \overline{M}_{\mathcal{E}}$ is the corresponding elementary embedding).
\end{prop}

\begin{proof}
See \cite{Kan:THI} 26.7.
\end{proof}

\section{The Product Reflection  Principle}

For any set $S$ of relational structures $\mathcal{A}=\langle A, \ldots \rangle$ of the same type,  the set-theoretic product $\prod S$  is the structure whose universe is the  set of all functions $f$ with domain $S$ such that $f(\mathcal{A})\in A$, for every $\mathcal{A}\in S$, and whose relations are defined pointwise.

\begin{definition}[The Product Reflection Principle ($\PRP$)]
For $\Gamma$  a definability class (i.e., one of $\Sigma_n$, $\Pi_n$, some $n>0$), and a set or a class $P$, $\Gamma(P)$-$\PRP$ asserts that for every $\Gamma$-definable, with parameters in $P$, proper class $\Ce$ of graphs the following holds:
\begin{enumerate}
\item[$\PRP$:] There is a subset $S$ of $\Ce$ such that for every  $G$ in $\Ce$ there is a homomorphism $\prod S \to G$.
\end{enumerate}
If $P=\emptyset$, then we simply write $\Gamma$-$\PRP$. If $P=V$, then we write $\Gamma$ in boldface, e.g., $\mathbf{\Sigma_n}$-$\PRP$.
\end{definition}

In the definition of $\Gamma (P)$-$\PRP$ we may replace ``graphs" by ``structures in a fixed (many-sorted, infinitary) relational language" and obtain equivalent principles, provided we allow for some additional parameters (see our remarks after definition \ref{defWVP}). Thus, the boldface principle $\mathbf{\Sigma_n}$-$\PRP$ for classes of graphs is equivalent to its version for classes of relational structures. 

\medskip

We shall denote by $C^{(n)}$  the $\Pi_n$-definable closed and unbounded class of ordinals $\kappa$ that are \emph{$\Sigma_n$-correct in $V$}, i.e., $V_\kappa\preceq_{\Sigma_n}V$. (See \cite{Ba:CC}.) 

\medskip

\begin{proposition}
\label{prop0}
$\mathbf{\Sigma_1}$-$\PRP$ holds. In fact, for every $\kappa \in C^{(1)}$ and every $\Sigma_1$-definable  with parameters in $V_\kappa$ proper class $\Ce$ of structures in a fixed  relational language $\tau \in V_\kappa$, the set $S:=\Ce\cap V_\kappa$ witnesses $\Sigma_1(V_\kappa)$-$\PRP$.
\end{proposition}

\begin{proof}
Let $\kappa \in C^{(1)}$, and let $\Ce$ be a $\Sigma_1$-definable, with a set of parameters $P\in V_\kappa$, proper class of structures in a relational language $\tau \in V_\kappa$. Note that since $\kappa \in C^{(1)}$, $V_\kappa =H_\kappa$, hence $|TC(\{ \tau\}\cup P)|<\kappa$. Let $\varphi (x)$ be a $\Sigma_1$ formula, with parameters in $P$, defining $\Ce$. We claim that $S:=\Ce\cap V_\kappa$ satisfies $\PRP$. Given $\mathcal{A}\in \Ce$, let $\lambda\in C^{(1)}$ be greater than $\kappa$ and such that $\mathcal{A}\in V_\lambda$. Let $N\preceq V_\lambda$ be of cardinality less than $\kappa$ and such that $\mathcal{A}\in N$ and $TC(\{ \tau\}\cup P)\subseteq N$. Let $\pi:M\to N$ be the inverse transitive collapse isomorphism, and let $\mathcal{B}\in M$ be such that $\pi(\mathcal{B})=\mathcal{A}$. Notice that $\pi$ fixes $\tau$ and the parameters of $\varphi(x)$. Since $M$ is transitive and of cardinality less than $\kappa$, $\mathcal{B}\in H_\kappa =V_\kappa$. Also, since $V_\lambda \models \varphi(\mathcal{A})$, we have  $N\models \varphi(\mathcal{A})$, and therefore $M\models \varphi(\mathcal{B})$. Hence, since $M$ is transitive and $\varphi$ is upwards absolute for transitive sets,  $\mathcal{B}\in \Ce$. Thus, $\mathcal{B}\in S$. Then the composition of $\pi$ with the projection $\prod S\to \mathcal{B}$ yields the desired homomorphism.
\end{proof}

\begin{proposition}
\label{prop1}
If $\kappa$ is a strong cardinal, then $\Sigma_2(V_\kappa)$-$\PRP$ holds.
\end{proposition}

\begin{proof}
Let $\kappa$ be a strong cardinal and let $\Ce$ be a $\Sigma_2$-definable, with parameters in $V_\kappa$,  proper class of structures in a fixed  relational language $\tau\in V_\kappa$. Let $\varphi (x)$ be a $\Sigma_2$ formula defining it. We will show that $S:=\Ce \cap V_\kappa$ witnesses $\PRP$. 

Given any $\mathcal{A}\in \Ce$, let $\lambda\in C^{(2)}$ be greater than or equal to $\kappa$ and with  $\mathcal{A}\in V_\lambda$.

Let $j:V\to M$ be an elementary embedding, with $crit(j)=\kappa$, $V_\lambda \subseteq M$, and $j(\kappa)>\lambda$.

By elementarity, the restriction of $j$ to $\prod S$ yields a homomorphism
$$h:\prod S\to \prod (\{ X: M\models \varphi(X)\} \cap V^M_{j(\kappa)}).$$
Since $\mathcal{A}\in V_\lambda$, and $\lambda \in C^{(2)}$, we have that $V_\lambda \models \varphi(\mathcal{A})$. Hence, since the fact that $\lambda \in C^{(1)}$ is $\Pi_1$-expressible and therefore downwards absolute for transitive classes, and since $V_\lambda \subseteq M$, it follows  that  $V_\lambda \preceq_{\Sigma_1}M$ and therefore  $M\models \varphi(\mathcal{A})$. Moreover $\mathcal{A}\in V_\lambda \subseteq V^M_{j(\kappa)}$. Thus, letting
$$g: \prod (\{ X: M\models \varphi(X)\} \cap V^M_{j(\kappa)})\to \mathcal{A}$$
be the projection map, we have that 
$$g\circ h:\prod S\to \mathcal{A}$$
is a homomorphism, as wanted.
\end{proof}

\begin{corollary}
\label{coro1}
If there exist a proper class of strong cardinals, then ${\mathbf \Sigma_2}$-$\PRP$ holds.
\end{corollary}

 \medskip
 
 We shall next show that $\SWVP$ is equivalent to the assertion that $\PRP$ holds for all definable proper classes of structures.
 
 \begin{proposition}
\label{prop3}
$\Gamma_n (P)$-$\PRP$ implies $\Gamma_{n}(P)$-$\SWVP$, for all $n>0$ and  every  $P$. 
\end{proposition}
 
\begin{proof}
Assume  $\mathcal{G}=\langle G_\alpha :\alpha \in \OR\rangle$ is a $\Gamma_n$-definable, with parameter $p\in P$,  sequence of graphs  such that for every $\alpha \leq \beta$ there exists some homomorphism $G_\beta \to G_\alpha$. Without loss of generality, the sequence is injective, i.e., $G_\alpha \ne G_\beta$ whenever $\alpha \ne \beta$. We shall produce a homomorphism  $G_\alpha \to G_\beta$, for some $\alpha <\beta$. Notice that the collection of all $G$ such that $G=G_\alpha$, some $\alpha \in \OR$, is a proper class. Let $\varphi(x)$ be a $\Gamma_n$ formula, with $p$ as a parameter, that defines $\mathcal{G}$.
 
 Let $\He$ be the class of all $G$ such that $G=G_\alpha$, for some ordinal $\alpha$ such that the rank of $G_\alpha$ is greater than or equal to $\alpha$. Notice that $\He$ is also a proper class. We claim that $\He$ is $\Gamma_n$-definable with $p$ as a parameter. This is clear if $\Gamma=\Sigma$, because $G\in \mathcal{H}$ iff there exists $\alpha$ such that $\varphi(\langle \alpha , G\rangle)$ and $\rm{rank}(G)\geq \alpha$. If $\Gamma=\Pi$, and $n=1$, then this is also clear since  $G\in \mathcal{H}$ iff for all transitive $M$ with $G, p\in M$, if $M\models$ ``$\{ x :\varphi(x)\}$ is a sequence of graphs indexed by the ordinals", then $M\models$ ``There exists $\alpha$ such that $\varphi(\langle \alpha , G\rangle)$ and $\rm{rank}(G)\geq \alpha$". If $\Gamma=\Pi$, and $n>1$, then $G\in \mathcal{H}$ iff for all $\beta \in C^{(n-1)}$ with $G, p \in V_\beta$, if $V_\beta \models$ `` $\{ x :\varphi(x)\}$ is a sequence of graphs indexed by the ordinals", then $V_\beta \models$ ``There exists $\alpha$ such that $\varphi(\langle \alpha , G\rangle)$ and $\rm{rank}(G)\geq \alpha$".
 
 Hence  by $\PRP$ there is a subset $S$ of $\He$  and, for every $G\in \He$, a homomorphism $\prod S\to G$.   Thus, for every $\beta$, there is a homomorphism $$h_\beta :\prod S\to G_{\beta}.$$

Let $\kappa$ be  the supremum of all $\alpha$ such that $G_\alpha\in S$.   
Pick any $\alpha >\kappa$. By our assumption, for every $\gamma \leq \kappa$ there is a homomorphism  $k_\gamma : G_\alpha \to G_\gamma$. Then the map $\ell_\alpha: G_\alpha \to \prod S$, given by $\ell_\alpha(x)=\{ \langle G_\gamma ,  k_\gamma (x)\rangle :G_\gamma \in S \}$ is a homomorphism. 

Now, given any $\alpha <\beta$ with $\alpha$ greater than $\kappa$, the composition  $h_\beta \circ \ell_\alpha :G_\alpha \to G_\beta$ is a homomorphism, as wanted.
\end{proof}

The converse is also true. Namely,

\begin{proposition}
\label{prop2}
$\Gamma_{n}(P)$-$\SWVP$ implies $\Gamma_n (P)$-$\PRP$, for every $n>1$, and every $P$. 
\end{proposition}

\begin{proof}
Let $n>0$ and fix a $\Gamma_n$-definable, with parameter $p$ in $P$, proper class  $\Ce$ of graphs and, aiming for a contradiction, suppose that $\PRP$ fails for $\Ce$. We build by induction on $\gamma$ a  sequence  $\overline{\Ce}=\langle \overline{\Ce}_\gamma :\gamma \in \OR\rangle$, where $\overline{\Ce}_\gamma =\prod (\Ce \cap V_{\alpha_\gamma} )$, some $\alpha_\gamma$, such that $\gamma \leq \eta$ implies $\alpha_\gamma \leq \alpha_\eta$, and such that there is no homomorphism $h:\overline{\Ce}_\gamma \to \overline{\Ce}_\eta$, whenever $\gamma <\eta$. Namely,  let $\alpha_0$ be the least ordinal such that $\Ce\cap V_{\alpha_0} \ne \emptyset$ and let $\overline{\Ce_0}= \prod (\Ce \cap V_{\alpha_0})$. Given $\overline{\Ce}_\delta$ and $\alpha_\delta$ for $\delta <\gamma$, let  $\beta$ be the least ordinal greater than ${\rm sup}\{ \alpha_\delta : \delta <\gamma\}$ such that for some $\mathcal{A}\in \Ce \cap V_{\beta}$ there is no homomorphism
 $$\prod (\Ce \cap V_{{\rm sup}\{\alpha_\delta :\delta <\gamma\}})  \to \mathcal{A}.$$
 Then let $\alpha_\gamma =\beta$ and define  $$\overline{\Ce}_\gamma =\prod (\Ce\cap V_{\alpha_\gamma}) \, .$$
 
Since $n>1$,   it is easily seen that $\overline{\Ce}$ is $\Gamma_n$-definable (with parameter $p$). For if $\Ce$ is $\Pi_n$-definable, then $X=\overline{\Ce}_\gamma$ if and only if 
 for every ordinal $\xi\in C^{(n-1)}$ with $p\in V_\xi$, , if $X\in V_\xi$, then $V_\xi\models ``X=\overline{\Ce}_\gamma"$. And if $\Ce$ is $\Sigma_n$-definable, then  $X=\overline{\Ce}_\gamma$ if and only if 
 for some ordinal $\xi\in C^{(n-1)}$ with $p$ and $X\in V_\xi$, $V_\xi\models ``X=\overline{\Ce}_\gamma"$.

By the construction of $\overline{\Ce}$, there is no homomorphism $\overline{\Ce_\gamma} \to \overline{\Ce_\delta}$ whenever $\gamma <\delta$. So, by $\Gamma_n(P)$-$\SWVP$ there are $\gamma <\delta$ such that there is no homomorphism $\overline{\Ce}_\delta \to \overline{\Ce}_\gamma$. But the projection is such a homomorphism.

\medskip

Suppose now $\Ce$ is $\Pi_1$-definable.  
 Let $\mathcal{D}$ be the   class of all  structures 
 $$\mathcal{D}_\gamma:=\langle \prod_{\alpha \in [\alpha^0 , \alpha_\gamma)} \!\!V_{\alpha +1}, \overline{\in}  ,  E_\gamma   \rangle$$
 where 
 \begin{enumerate}
 \item $\overline{\in}$ is the pointwise membership relation, 
 \item 
 $$E_\gamma \subseteq \prod_{\alpha \in [\alpha^0 , \alpha_\gamma)}\!\! (\prod( \Ce \cap V_\alpha))$$
consists of all $f$ such that $f(\alpha)\subseteq f(\alpha')$, for all $\alpha <\alpha'$ in their domain, and
\item $\alpha^0$ is the least ordinal in $C^{(1)}$ such that $\Ce \cap V_{\alpha^0} \ne \emptyset$.
\item $\alpha_\gamma$ is the $\gamma$-th element of the class of ordinals $\eta$ in $C^{(1)}$ greater than $\alpha^0$ such that
$$\forall \delta <\eta  \exists    A\in \Ce \cap V_{\eta}  ( \neg \exists h  (h:\prod (\Ce \cap V_\delta) \to A \mbox{ is a homomorphism})).$$

\end{enumerate}

Since $\PRP$ fails for $\Ce$, $\De_\gamma$ exists for every ordinal $\gamma$. Also $\mathcal{D}$ is $\Pi_1$, for $X\in \mathcal{D}$ if and only if every transitive model  of a sufficiently rich finite fragment of ZFC that contains $X$ satisfies ``$X\in \mathcal{D}$".

If $\gamma <\gamma'$, then there is no homomorphism $h:\mathcal{D}_\gamma \to \mathcal{D}_{\gamma'}$, for the restriction  of such an $h$ to $E$  would yield a homomorphism $E_\gamma \to E_{\gamma'}$, 
which in turn would  yield  (by (2) above)  a homomorphism
$$h': \prod (\Ce \cap V_{\alpha_{\gamma}}) \to \prod (\Ce \cap V_{\alpha_{\gamma'}})$$
by letting $$h'(f)= \bigcup (h(\langle f\restriction V_\alpha\rangle_{\alpha \in [\alpha^0 , \alpha_\gamma)})).$$ 
But by (3)-(b) above, there is some $A\in \Ce \cap V_{\alpha_{\gamma'}}$ for which there is no homomorphism
$$\prod (\Ce \cap V_{\alpha_\gamma})\to A.$$

So, by $\Pi_1(P)$-$\SWVP$ there must exist $\gamma <\gamma'$ for which there is no homomorphism $\mathcal{D}_{\gamma'}\to \mathcal{D}_\gamma$. But the projection is such a homomorphism.
 \end{proof}

The following is an immediate corollary of Propositions \ref{prop2} and \ref{prop3}.
\begin{corollary}
\label{coro2}
$\Gamma_n (P)$-$\PRP$ is equivalent to $\Gamma_{n}(P)$-$\SWVP$, for every $n>0$.  Hence, $\SWVP$ is equivalent to $\mathbf{\Sigma_n}$-$\PRP$ holds for every $n$.
\end{corollary}

\section{The equivalence of $\WVP$ and $\SWVP$}

Following the proof due to the second author \cite{W:WVP} of the equivalence of $\WVP$ and $\SWVP$, we prove next the equivalence of their corresponding versions for definable classes of graphs. We only give a sketch of the proof, emphasizing only the definability aspects of it, and we refer to \cite{W:WVP} for additional details.

\begin{theorem}
\label{equivWandSW}
For $\Gamma$ a definability class, and every set or proper class $P$, the principles $\Gamma(P)$-$\WVP$ and $\Gamma(P)$-$\SWVP$ are equivalent.
\end{theorem}

\begin{proof}
That $\Gamma(P)$-$\SWVP$ implies $\Gamma(P)$-$\WVP$ is immediate. For the converse, fix $\Gamma$ and $P$, and assume that $\langle G_\alpha :\alpha \in \OR\rangle$ is a sequence of graphs that is a counterexample to $\Gamma(P)$-$\SWVP$. Thus, the sequence is $\Gamma(P)$-definable,  for every $\alpha \leq \alpha'$ there is a homomorphism from $G_{\alpha'}$ to $G_\alpha$, and for every $\alpha <\alpha'$ there is no homorphism from $G_\alpha$ to $G_{\alpha'}$.  

Let 
$$
    H_\beta =
    \begin{cases*}
      \prod_{\alpha \leq \beta}G_\alpha & if $\beta=0$ or $\beta$ is a successor ordinal,\\
     \prod_{\alpha <\beta}G_\alpha       & if $\beta$ is a limit ordinal.
    \end{cases*}
 $$
  We claim that the sequence $\langle H_{\beta} :\beta \in \OR\rangle$ yields also a counterexample to $\Gamma(P)$-$\SWVP$. First, since $\Sigma_1(P)$-$\SWVP$ is provable in ZFC (see propositions \ref{prop0}, \ref{prop3}), we may assume $\Gamma$ is either $\Sigma_n$ with $n>1$, or $\Pi_n$ with $n\geq 1$. In either case, since $\langle G_\alpha :\alpha \in \OR\rangle$ is $\Gamma(P)$-definable, so is $\langle H_{\beta} :\beta \in \OR\rangle$ (see the proof of \ref{prop2}). Second, for every $\beta \leq \beta'$ there is a restriction homomorphism from $H_{\beta'}$ to $H_\beta$, and for every $\beta <\beta'$ there is no homomorphism from $H_\beta$ to $H_{\beta'}$, for otherwise, arguing similarly as in proposition \ref{prop3}, we could compose homomorphisms
$$G_\beta \to H_\beta \to H_{\beta'}\to G_{\beta +1}$$
to obtain a homomorphism from $G_\beta$ to $G_{\beta +1}$, thus yielding a contradiction.

The class $\Lambda$ of limit ordinals $\lambda$ such that $G_\alpha \in V_\lambda$ for all $\alpha <\lambda$ is closed and unbounded. If the sequence $\langle G_\alpha :\alpha \in \OR\rangle$  is $\Pi_1(P)$-definable, then so is $\Lambda$, and if it is $\Sigma_n(P)$-definable, with $n>1$, then so is $\Lambda$. Now for each $\lambda \leq \lambda'$ in $\Lambda$, define the function $h_{\lambda' ,\lambda}:V_{\lambda' +1}\to V_{\lambda +1}$ by
$$h_{\lambda ,\lambda'}(x)=x\cap V_\lambda.$$
For each $\lambda \in \Lambda$ let the structure 
$$\mathcal{M}_\lambda =\langle V_{\lambda +1},\in, \lambda, R_\lambda, S_\lambda, T_\lambda\rangle$$
be such that 
\begin{enumerate}
\item[] $R_\lambda =\{ \langle \beta , x,y\rangle: (\beta =$ rank$(x)$ and $x\in y)$ or $\beta =\lambda\}$
\item[] $S_\lambda =\{ \langle \beta , x,y\rangle: (\beta =$ rank$(x)$ and $x\not \in y)$ or $\beta =\lambda\}$
\item[] $R_\lambda =\{ \langle \beta , x,y\rangle: x$ is adjacent to $y$ in $H_\beta \}$
\end{enumerate}
If the sequence $\langle G_\alpha :\alpha \in \OR\rangle$  is $\Pi_1(P)$-definable, then as argued above so are the sequence $\langle H_{\beta} :\beta \in \OR\rangle$ and  $\Lambda$, hence also $\langle \mathcal{M}_\lambda :\lambda \in \Lambda\rangle$; and if it is $\Sigma_n(P)$-definable, with $n>1$, then so is $\langle \mathcal{M}_\lambda :\lambda \in \Lambda\rangle$. 

Now as in \cite{W:WVP} one can show that for all $\lambda \leq \lambda'$ in $\Lambda$, the map $h_{\lambda',\lambda}:\mathcal{M}_{\lambda'}\to \mathcal{M}_\lambda$ is a homomorphism, and is the unique one. Moreover, for every $\lambda <\lambda'$ in $\Lambda$ there is no homomorphism from $\mathcal{M}_\lambda$ to $\mathcal{M}_{\lambda'}$. 
By re-enumerating the class $\Lambda$  in increasing order as $\langle \lambda_\alpha :\alpha \in \OR\rangle$ we obtain  a sequence of relational structures $\langle \mathcal{M}_{\lambda_\alpha}:\alpha \in \OR\rangle$ which yields a counterexample to $\Gamma(P)$-$\WVP$.
\end{proof}

\section{The main theorem for strong cardinals}

\begin{theorem}
\label{thm1}
The following are equivalent:
\begin{enumerate}
\item There exists a  strong cardinal.
\item $\Sigma_2$-$\PRP$
\item $\Pi_1$-$\PRP$
\item $\Sigma_2$-$\SWVP$
\item $\Pi_1$-$\SWVP$
\item $\Sigma_2$-$\WVP$
\item $\Pi_1$-$\WVP$.
\end{enumerate}
\end{theorem}

\begin{proof}
(1)$\Rightarrow$(2) is given by proposition \ref{prop1}; (2)$\Rightarrow$(3), (4)$\Rightarrow$(5), and (6)$\Rightarrow$(7) are   immediate; the equivalence of (2) and (4), and also of (3) and (5), are given by corollary \ref{coro2}. The equivalence of (4) and (6), and also of (5) and (7), is given by theorem \ref{equivWandSW}.  So, it will be sufficient to prove (3)$\Rightarrow$(1).

\medskip

(3)$\Rightarrow$(1):  
Let $\mathcal{A}$ be the class of all structures  
$$\Ae_\alpha:=\langle V_{ \alpha+1}, \in ,  \alpha  ,\{ R^{\alpha}_\varphi \}_{\varphi \in \Pi_1} \rangle$$   where  the constant $ \alpha$ is the $\alpha$-th element of $C^{(1)}$ 
 and   $\{ R^{\alpha}_\varphi\}_{\varphi \in \Pi_1}$ is the $\Pi_1$ relational diagram for $V_{\alpha +1}$, i.e., if $\varphi(x_1,\ldots ,x_n)$ is a $\Pi_1$ formula in the language of $\langle V_{ \alpha +1},\in,\alpha   \rangle$, then $$R^\alpha_\varphi =\{ \langle x_1,\ldots ,x_n\rangle: \langle V_{ \alpha +1},\in,\alpha    \rangle\models ``\varphi(x_1,\ldots ,x_n)"\}\, .$$

\medskip

We claim that $\mathcal{A}$ is $\Pi_1$-definable without parameters. For $X\in \mathcal{A}$ if and only if $X=\langle X_0,X_1,X_2,X_3  \rangle$, where

 \begin{enumerate}
 \item[(1)] $X_2$  belongs to $C^{(1)}$ 
 \item[(2)] $X_0=V_{X_2 +1}$
 \item[(3)] $X_1 =\in \restriction X_0$
 \item[(4)] $X_3$ is the $\Pi_1$ relational diagram of $\langle X_0,X_1,X_2  \rangle$, and 
 \item[(5)] $\langle X_0, X_1, X_2 \rangle\models ``X_2$ is the $X_2$-th element of $C^{(1)}$".
   \end{enumerate}

\medskip
    
Note that $\mathcal{A}$ is a proper class. In fact, the class $C$ of ordinals $\alpha$ such that $\Ae_\alpha \in \Ae$ is a closed and unbounded proper class..
By $\Pi_1$-$\PRP$  there exists a subset $S$ of $C$ such that for every $\beta\in C$ there is a homomorphism   $j_\beta : \prod_{\alpha \in S} \Ae_\alpha\to \Ae_\beta$. By enlarging $S$, if necessary, we may assume that $\rm{sup}(S)\in S$.
Let us denote $\prod_{\alpha \in S}\Ae_\alpha$ by $M$. Notice that   
$$M =\langle \prod_{\alpha \in S} V_{\alpha +1}, \overline{\in}, \langle \alpha \rangle_{\alpha \in S},  \{ \overline{R}^\alpha_\varphi\}_{\varphi \in \Pi_1} \rangle$$
where $\overline{\in}$ is the pointwise membership relation, and $\overline{R}^\alpha_\varphi$ is the pointwise $R^\alpha_\varphi$ relation. Let $\kappa:= \rm{sup}(S)$.  

\medskip

Now fix some $\beta\in C$  greater than $\kappa$,  of uncountable cofinality, and assume, aiming for a contradiction, that no cardinal $\leq \kappa$ is $\beta$-strong. Let $j=j_\beta$. 

\begin{claim}
$j$ preserves the Boolean operations $\cap$, $\cup$, $-$, and also the $\subseteq$ relation.
\end{claim}

\begin{proof}[Proof of claim]
For every $X,Y,Z\in M$, $$M\models ``X=Y\cap Z"\quad \mbox{ iff }\quad V_{\alpha +1}\models ``X(\alpha)=Y(\alpha)\cap Z(\alpha)",   \mbox{ all } \alpha \in S\, .$$ So, letting $\varphi (x,y,z)$ be the bounded formula  expressing   $x=y\cap z$, we have that $\langle X(\alpha),Y(\alpha),Z(\alpha)\rangle \in R^\alpha_\varphi$, for all $\alpha \in S$. Hence $\langle X,Y,Z\rangle \in \overline{R}_\varphi$, and since $j$ is a homomorphism $\langle j(X),j(Y),j(Z)\rangle \in R^\beta_\varphi$, which yields  $\Ae_\beta \models ``j(X)=j(Y)\cap j(Z)"$. 

Similarly for the operations $\cup, -$, and for the relation $\subseteq$. 
\end{proof}

Now define
$k:V_{\kappa +1}\to V_{\beta +1}$
by
$$k(X)=j(\langle X\cap V_{\alpha}\rangle_{\alpha \in S})\, .$$

\begin{claim}
$k$ also preserves the Boolean operations, as well as the $\subseteq$ relation. 
\end{claim}

\begin{proof}[Proof of claim]
Suppose $V_{\kappa +1}\models ``X=Y\cap Z"$. Then $X\cap V_{\alpha} =(Y\cap V_{\alpha})\cap (Z\cap V_{\alpha})$, for every $\alpha \in S$. Hence, 
$$M\models ``\langle X\cap V_{\alpha} \rangle_{\alpha \in S} =\langle Y\cap V_{\alpha}\rangle_{\alpha \in S}\cap \langle Z\cap V_{\alpha}\rangle_{\alpha \in S}"\, .$$
Since $j$ preserves the $\cap$ operation,
$$\Ae_{\beta}\models ``k(X)=k(Y)\cap k(Z)"\, .$$
Let $\psi(x)$ be the $\Pi_1$ formula (in the language with an additional constant for $\alpha$) asserting that $x=V_{\alpha}$. Then 
$$\langle  V_{\alpha} \rangle_{\alpha \in S}\in \overline{R}_\psi $$ and so $j(\langle  V_{\alpha} \rangle_{\alpha \in S})\in R^\beta_\psi$, which yields $j(\langle  V_{\alpha} \rangle_{\alpha \in S})=V_{\beta}$. Hence, since $j$ preserves the subset relation, and sends $\langle \alpha\rangle_{\alpha \in S}$ to $\beta$,   we have  that $k(X), k(Y), k(Z) \subseteq V_\beta$. So,  
$$V_{\beta +1}\models ``k(X)=k(Y)\cap k(Z)"\, .$$
Similarly for the operations $\cup , -$, and the relation $\subseteq$.
\end{proof}

\begin{claim}
$k$ maps ordinals to ordinals, and is the identity on $\omega+1$.
\end{claim}

\begin{proof}[Proof of claim]
Let $\varphi(x)$ be the bounded formula expressing that $x$ is an ordinal. Let $\gamma \leq \kappa$. Then $\gamma \cap V_{\alpha}$ is an ordinal, for all $\alpha <\kappa$, and so  
$$M\models ``\langle \gamma \cap V_{\alpha}\rangle_{\alpha \in S} \, \overline{\in} \, \, \overline{R}_\varphi"\, .$$
Since $j$ is a homomorphism,
$$\Ae_\beta \models ``j(\langle \gamma \cap V_{\alpha}\rangle_{\alpha \in S}) \in R^\beta_\varphi"$$
which yields that $k(\gamma)=j(\langle \gamma \cap V_{\alpha}\rangle_{\alpha \in S})$ is an ordinal in $\Ae_\beta$, hence also in $V_{\beta +1}$. 

For every ordinal $\gamma \leq \omega$, we have that $\gamma \cap V_{\alpha} =\gamma$, for all $\alpha \in S$. Moreover, $\gamma$ is definable by some bounded formula $\varphi_\gamma$. Hence,
$$M\models ``\langle \gamma \cap V_{\alpha}\rangle_{\alpha \in S} \, \overline{\in}\, \,  \overline{R}_{\varphi_\gamma}"$$
and therefore
$$\Ae_\beta \models ``j(\langle \gamma \cap V_{\alpha}\rangle_{\alpha \in S}) \in R^\beta_{\varphi_\gamma}"$$
which yields
$k(\gamma)=\gamma$.
\end{proof}

\medskip

Note that $k(\kappa)=j(\langle \alpha\rangle_{\alpha \in S})= \beta$.

\medskip

For each $a\in [\beta]^{<\omega}$,   define $E_a$ by
$$X\in E_a \quad \mbox{ iff }\quad X\subseteq [\kappa]^{|a|}  \mbox{ and } a\in k(X)\, .$$
Since $k(\kappa)=\beta$ and $k(|a|)=|a|$, we also have $k([\kappa]^{|a|})=[\beta]^{|a|}$, hence $[\kappa]^{|a|} \in E_a$. Moreover, since $k$ preserves Boolean operations and the $\subseteq$ relation, $E_a$ is an ultrafilter over $[\kappa]^{|a|}$. 

\begin{claim}
$E_a$ is  $\omega_1$-complete.
\end{claim}

\begin{proof}[Proof of claim]
Given $\{ X_n:n<\omega \}\subseteq E_a$, let $Y=\{ \langle n,x\rangle: x\in X_n\}$. So, $Y\subseteq V_\kappa$. We can express that  $X=\bigcap_{n<\omega}X_n$ by a bounded sentence $\varphi$ in the  parameters  $X$, $Y$ and $\omega$. Moreover, since $\alpha$ is a limit ordinal, for every $\alpha \in S$, the sentence $\varphi(X\cap V_{\alpha}, Y\cap V_{\alpha}, \omega)$ expresses that $X\cap V_{\alpha} =\bigcap_{n<\omega}X_n\cap V_{\alpha}$.  So, 
$$M\models ``\langle  X \cap V_{\alpha}, Y\cap V_{\alpha},  \omega \rangle_{\alpha \in S} \, \overline{\in}\, \,  \overline{R}_\varphi"\, .$$
Since $j$ is a homomorphism,
$$\Ae_\beta \models ``\langle j( X \cap V_{\alpha}), j(Y\cap V_{\alpha}), j(\omega) \rangle_{\alpha \in S} \in R^\beta_\varphi"$$
and so $\langle k(X),k(Y),k(\omega)\rangle$ satisfies $\varphi$. Since $k(\omega)=\omega$, we thus have $k(X)=\bigcap_{n<\omega}k(X_n)$. Hence, $a\in k(X)$, and so $X\in E_a$.
\end{proof} 

Let $\mathcal{E}:=\{ E_a: a\in [\beta]^{<\omega}\}$.

\begin{claim}
\label{normal}
$\mathcal{E}$ is normal. That is, whenever $a\in [\beta]^{<\omega}$ and $f$ is a function with domain $[\kappa]^{|a|}$  such that $\{ s\in [\kappa]^{|a|}: f(s)\in {\rm{max}}(s)\}\in E_a$, there is $b\supseteq a$   such that
$\{ s\in [\kappa]^{|b|}: f(\pi^\kappa_{ba}(s))\in s\}\in E_{b}$, where $\pi^\kappa_{ba}:[\kappa]^{|b|}\to [\kappa]^{|a|}$ is the standard projection function.
\end{claim}

\begin{proof}
Fix $a$ and $f$, and suppose 
$$X:= \{ s\in [\kappa]^{|a|}: f(s)\in {\rm max}(s)\}\in E_a \, .$$
Since the formula $\varphi(X, \kappa, |a|, f)$ defining $X$ is  bounded, and since $k(\kappa)=\beta$ and $k(|a|)=|a|$, we have that
$$k(X)=\{ s\in[\beta]^{|a|}: k(f)(s)\in {\rm max}(s)\}\, .$$
Also, since $X\in E_a$, we have that $k(f)(a)\in {\rm max}(a)$. 

Let $\delta =k(f)(a)$, and let $b=a\cup \{ \delta \}$. Thus,
$$b\in \{ s\in [\beta]^{|b|}: k(f)(\pi^\beta_{ba}(s))\in s\}$$
where $\pi^\beta_{ba}:[\beta]^{|b|}\to [\beta]^{|a|}$ is the standard projection function.
So, since $\{ s\in [\beta]^{|b|}: k(f)(\pi^\beta_{ba}(s))\in s\}=k(\{ s\in [\kappa]^{|b|}: f(\pi^\kappa_{ba}(s))\in s\})$, we have 
$$\{ s\in [\kappa]^{|b|}: f(\pi^\kappa_{ba}(s))\in s\}\in E_b$$
which shows that $\Ee$ is normal.
\end{proof}

For each $a\in [\beta]^{<\omega}$, the ultrapower ${\rm{Ult}}(V,E_a)$ of $V$ by the $\omega_1$-complete ultrafilter $E_a$ is well-founded. So, let 
$$j_a:V\to M_a \cong {\rm{Ult}}(V, E_a)$$
with $M_a$ transitive, be the corresponding ultrapower embedding. As usual, we denote the elements of $M_a$ by their corresponding elements in ${\rm Ult}(V, E_a)$.

\begin{claim}
$\mathcal{E}$ is coherent. I.e., for every $a\subseteq b$ in $[\beta]^{<\omega}$,
$$X\in E_a \mbox{ if and only if } \{ s\in [\kappa]^{|b|}: \pi_{ba}(s) \in X\}\in E_b\, .$$
\end{claim}

\begin{proof}
Let $a\subseteq b$ in $[\beta]^{<\omega}$, and suppose $X\in E_a$. Thus, $X\subseteq [\kappa]^{|a|}$ and $a\in k(X)$. We need to see that $b\in k(\{ s\in [\kappa]^{|b|}:\pi_{ba}(s) \in X\})$. Now notice that, since $k$ is the identity on natural numbers, and $k(\kappa)=\beta$,  $$k(\{ s\in [\kappa]^{|b|}:\pi_{ba}(s) \in X\})=\{ s\in [\beta]^{|b|}: \pi_{ba}(s)\in k(X)\}.$$
Hence, since $\pi_{ba}(b)=a$, and $a\in k(X)$, we have that $b\in \{ s\in [\beta]^{|b|}: \pi_{ba}(s)\in k(X)\}$, as wanted.

Conversely, if $\{ s\in [\kappa]^{|b|}: \pi_{ba}(s) \in X\}\in E_b$, we have that $b\in k(\{ s\in [\kappa]^{|b|}:\pi_{ba}(s) \in X\})=\{ s\in [\beta]^{|b|}: \pi_{ba}(s)\in k(X)\}.$ Hence, $\pi_{ba}(b)=a\in k(X)$, and therefore  $X\in E_a$.
\end{proof}

\medskip

For each $a\subseteq b$ in $[\beta]^{<\omega}$, let $i_{ab}:M_a\to M_b$ be  given by 
$$i_{ab}([f]_{E_a})=[f\circ \pi_{ba}]_{E_b}$$
for all $f:[\kappa]^{|a|}\to V$. By coherence, the maps $i_{ab}$ are  well-defined and commute with the ultrapower embeddings $j_a$ (see \cite{Kan:THI} 26). 

Let   $M_\Ee$ be  the  direct limit of $$\langle \langle M_a:a\in [\beta]^{<\omega}\rangle, \langle i_{ab}:a\subseteq b\rangle\rangle.$$   

For notational simplicity, whenever we write $[a, [f]]\in [b,[g]]$ in $M_\Ee$, what we mean is that $[f]=[f]_{E_a}\in M_a$, $[g]=[g]_{E_b}\in M_b$, and $[\langle a, [f]_{E_a}\rangle ]]_\Ee \in_\Ee [\langle b, [g]_{E_b}\rangle ]_\Ee$.

Let $j_\Ee :V\to M_\Ee$ be the corresponding limit elementary embedding, i.e., 
$$j_\Ee(x)=[a,[c^a_x]_{E_a}]$$
for some (any) $a\in [\beta]^{<\omega}$, and  where  $c^a_x:[\kappa]^{|a|}\to \{ x\}$. 

Let   $k_a:M_a\to M_\Ee$ be given by $$k_a([f]_{E_a})=[a, [f]_{E_a}] \, .$$
It is easily checked that $j_\Ee =k_a \circ j_a$ and $k_b \circ i_{ab}=k_a$, for all $a\subseteq b$, $a,b \in [\beta]^{<\omega}$. 
Thus,  letting ${\rm Id}_{|a|}: [\kappa]^{|a|}\to [\kappa]^{|a|}$ be the identity function, we have 
$$M_\Ee =\{ j_\Ee (f)(k_a([{\rm id}_{|a|}]_{E_a})) :  a\in [\beta]^{<\omega} \mbox{  and  } f  :  [\kappa]^{|a|}\to V\}\, .$$

\medskip

Let $M^\ast_\Ee  :=\{ [a,[f]]\in M_\Ee : f \cap V_{\alpha}: [\alpha]^{|a|}\to V_{\alpha} , \mbox{ all }\alpha \in S\}$.
Suppose $[a,[f]], [b,[g]]\in M^\ast_\Ee$. Then the following can be easily verified:
\begin{enumerate}
\item $[a,[f]]\in_\Ee  [b,[g]] \mbox{  iff  } k(f)(a)\in k(g)(b)$
\item $[a,[f]] =_\Ee  [b,[g]] \mbox{  iff  } k(f)(a) = k(g)(b)$
\end{enumerate}

\begin{claim}
\label{wellfounded}
$M^\ast_\Ee$ is well-founded and closed under  $\in_\Ee$.
\end{claim}

\begin{proof}
Well-foundedness follows from items (1) and (2) above, as any infinite $\in_\Ee$-descending sequence in $M^\ast_\Ee$ would yield an infinite $\in$-des\-cend\-ing sequence in $V_{\beta +1}$.

Now suppose $[a,[f]]\in_\Ee [b,[g]]$, with $[b,[g]]\in M^\ast_\Ee$. Then for some $c\supseteq a,b$, and some $X\in E_c$,
$$(f\circ \pi_{ca})(s)\in (g\circ \pi_{cb})(s)$$
for every $s\in X$. Let $Y=\{ \pi_{ca}(s):s\in X\} \in E_a$. Define $h:[\kappa]^{|a|}\to V$ by:  $h(s)=f(s)$ for all $s\in Y$, and $h(s)=0$, otherwise. Then $[h]_{E_a}=[f]_{E_a}$, and $[a,[f]]=[a,[h]]\in M^\ast_\Ee$.
\end{proof}

By the last Claim, $M^\ast_\Ee$ is well-founded and extensional. So, let $M^\ast$ be the transitive collapse of  $M^\ast_\Ee$.

\begin{claim}
\label{betastrong}
$V_\beta \subseteq M^\ast$.
\end{claim}

\begin{proof}[Proof of claim]
Since $\kappa$ and $\alpha$, for $\alpha \in S$, belong to $C^{(1)}$, we have that  $|V_\kappa |=\kappa$ and $|V_{\alpha} |=\alpha$, all $\alpha \in S$. Let $f\in V$ be  a bijection between $[\kappa]^1$ and $V_\kappa$  such that $f\restriction [\alpha]^1$ is a bijection between $[\alpha]^1$ and $V_{\alpha}$, all $\alpha \in S$. 
Let $\varphi (x,y,z)$ be a $\Pi_1$ formula expressing that $x=[u]^1$, with $u$  an ordinal, $y=V_u$, and $z:x\to V_u$ is a bijection. Thus,
$$M\models ``\langle [\alpha]^1 , V_{\alpha}, f\cap V_{\alpha}\rangle_{\alpha \in S} \, \overline{\in}\, \, \overline{R}_\varphi"\, .$$
Hence,
$$\Ae_\beta \models ``\langle k([\kappa]^1), k(V_\kappa), k(f)\rangle \in R^\beta_\varphi"$$
and so $k(f): [\beta]^1 \to V_\beta$ is a bijection. Therefore, for every $x\in V_\beta$ there exists $\gamma <\beta$ such that $k(f)(\{ \gamma\})=x$.

Thus, letting $D:=\{ [\{\gamma\}, [f]]:\gamma <\beta\}$, we have just shown that the  map $i:\langle D, \in_\Ee \restriction D\rangle  \to \langle V_\beta, \in \rangle$ given by
$$i([\{ \gamma\},[f]])=k(f)(\{\gamma\})$$
is onto. Moreover, if $[\{ \gamma\},[f]]\in_\Ee [\{ \delta \},[f]]$, then for some $X\in E_{\{\gamma,\delta\} }$, we have 
$$(f\circ \pi_{\{ \gamma,\delta\}\{\gamma\}})(s)\in (f\circ  \pi_{\{\gamma ,\delta\}\{\delta\}})(s)$$
for every $s\in X$. Letting $\varphi$ be the bounded formula expressing this, and  since $V_{\alpha}$ is closed under $f$, for every $\alpha$, we have
$$M\models ``\langle X\cap V_{\alpha}, (f\circ \pi_{\{ \gamma,\delta\}\{\gamma\}})\cap V_{\alpha}, (f\circ \pi_{\{\gamma ,\delta\}\{\delta\}})\cap V_{\alpha}\rangle \, \overline{\in}\,\, \overline{R}_\varphi "\, .$$
Hence, in $\Ae_{\beta}$, for every $s\in k(X)$,
$$(k(f)\circ \pi_{\{ \gamma,\delta\}\{\gamma\}})(s)\in (k(f)\circ  \pi_{\{\gamma ,\delta\}\{\delta\}})(s)\, .$$
In particular, since $\{\gamma,\delta\}\in k(X)$, 
$$k(f)(\{\gamma\})\in k(f)(\{\delta\})\, .$$

A similar argument shows that $i$ is one-to-one. Hence, $i$ is an isomorphism, and so $i$ is just the transitive collapsing map. Since $D\subseteq M^\ast_\Ee$, to conclude  that $V_\beta \subseteq M^\ast$ it will be sufficient to show that the transitive collapse of $D$ is the same as the restriction to $D$ of the transitive collapse of $M^*_E$. For this, it suffices to see that every $\in_E$-element of an element of $D$ is $=_E$-equal to an element of $D$. 
So, suppose $[\{ \gamma \}, [f]]\in D$ and $[a,[g]] \in_E [\{\gamma\}, [f]]$, with $[a,[g]] \in M^*_E$. Then $k(g)(a)\in k(f)({\gamma})$, by (1) above (just before Claim \ref{wellfounded}). Now $k(f):[\beta]^1\to V_\beta$
is surjective and $V_\beta$ is transitive, so there is some $\delta < \beta$ such that
$k(f)({\delta})=k(g)(a)$. Hence, by (2) above, $[\{\delta\},[f]] =_E [a,[g]]$. 
\end{proof}

\begin{claim}
$M_\Ee$ is closed under $\omega$-sequences, hence it is well-founded.
\end{claim}

\begin{proof}[Proof of claim] Let $\langle  j_\Ee (f_n)(k_{a_n}([{\rm id}_{|a_n|}]_{E_{a_n}}))\rangle_{n<\omega}$ be a sequence of elements of $M_\Ee$. On the one hand, the sequence $\langle j_\Ee (f_n)\rangle_{n<\omega} =j_\Ee(\langle f_n\rangle_{n<\omega})$ belongs to $M_\Ee$. On the other hand,   $k_{a_n}([{\rm Id}_{|a_n|}]_{E_{a_n}})=[a_n, [{\rm Id}_{|a_n|}]_{E_{a_n}}]$ belongs to $M^\ast_\Ee$, all $n<\omega$. Since  $\Ee$ is  normal (Claim \ref{normal}), as in \cite{Kan:THI} 26.2 (a) we can show that the transitive collapse of $[a_n, [{\rm Id}_{|a_n|}]_{E_{a_n}}]$ is precisely $a_n$. The sequence $\langle a_n\rangle_{n<\omega}$ belongs to $V_\beta$, because $\beta$ has uncountable cofinality. Hence,  since $V_\beta \subseteq M^\ast$, the preimage of  $\langle a_n\rangle_{n<\omega}\in V_\beta$ under the transitive collapsing map of $M^\ast_\Ee$ to $M^\ast$,  is precisely  the sequence $\langle k_{a_n}([{\rm Id}_{|a_n|}]_{E_{a_n}})\rangle_{n<\omega}$ and belongs to $M_\Ee$. It now follows that the sequence $\langle  j_\Ee (f_n)(k_{a_n}([{\rm id}_{|a_n|}]_{E_{a_n}}))\rangle_{n<\omega}$ is also in $M_\Ee$.
\end{proof}

Let $\pi:M_\Ee \to N$ be the transitive collapsing isomorphism, and let $j_N:V\to N$ be the corresponding elementary embedding, i.e., $j_N=\pi \circ j_\Ee$.

\begin{claim}
$j_N(\kappa)\geq \beta$.
\end{claim}

\begin{proof}[Proof of claim] Let $\alpha <\beta$. Let ${\rm Id}_1$ be the identity function on $[\kappa]^1$, and let $c_{[\kappa]^1}:[\kappa]^1 \to \{ [\kappa]^1\}$ and $c_\kappa :[\kappa]^1\to \{ \kappa \}$. In $M_{\{ \alpha \}}$, we have
$$[{\rm Id}_1]_{E_{\{ \alpha \}}}\in [c_{[\kappa]^1}]_{E_{\{ \alpha \}}}=[[c_\kappa]_{E_{\{ \alpha \}}}]^1=[j_{\{ \alpha \}}(\kappa)]^1$$
hence in $M_\Ee$,
$$k_{\{ \alpha \}}([{\rm Id}_1]_{E_{\{ \alpha \}}})\in k_{\{ \alpha \}}([j_{\{ \alpha \}}(\kappa)]^1) =[j_\Ee (\kappa)]^1$$
and therefore, since $\pi (k_{\{ \alpha \}}([{\rm Id}_1]_{E_{\{ \alpha \}}}))=\{ \alpha \}$, in $N$ we have
$$\{ \alpha \} \in [j_N(\kappa)]^1$$
that is, $\alpha <j_N(\kappa)$.
\end{proof}

Since $\beta >\kappa$, the last claim implies that the critical point of $j_N$ is less than or equal to $\kappa$.
Thus,  since $V_\beta \subseteq N$ (by Claim \ref{betastrong}), $j_N$ witnesses that the critical point of $j_N$ is a $\beta$-strong cardinal. But this  is in contradiction to  our choice of $\beta$. This completes the proof of theorem \ref{thm1}.
\end{proof}

The boldface version of theorem \ref{thm1}, i.e., with parameters, also holds by essentially the same arguments. Namely,

\begin{theorem}
\label{thm3}
The following are equivalent:
\begin{enumerate}
\item There exists a  proper class of strong cardinals.
\item $\mathbf{\Sigma_2}$-$\PRP$
\item $\mathbf{\Pi_1}$-$\PRP$
\item $\mathbf{\Sigma_2}$-$\SWVP$
\item $\mathbf{\Pi_1}$-$\SWVP$
\item $\mathbf{\Sigma_2}$-$\WVP$
\item $\mathbf{\Pi_1}$-$\WVP$

\end{enumerate}
\end{theorem}

For  the proof of (3) implies (1), in order to show that there exists a strong cardinal greater than a fixed ordinal $\gamma$ we need to  consider the  class of structures 
$$\Ae_\alpha:=\langle V_{ \alpha+1}, \in ,  \alpha  ,\{ R^{\alpha}_\varphi \}_{\varphi \in \Pi_1}, \langle \delta \rangle_{\delta <\gamma} \rangle$$
where $$\Ae_\alpha:=\langle V_{ \alpha+1}, \in ,  \alpha  ,\{ R^{\alpha}_\varphi \}_{\varphi \in \Pi_1} \rangle$$ is as in the proof of Theorem \ref{thm1}, and we have a constant $\delta$ for every $\delta <\gamma$.

\section{The general case}

We shall now consider the general case of definable proper classes of structures with any degree of definable complexity. For this we shall need  the following new kind of large cardinals.

If $j:V\to M$ is an elementary embedding, with $M$ transitive and critical point $\kappa$, and $A$ is a class  definable by a formula $\varphi$ (possibly with parameters in $V_\kappa$), we define
$$j(A):=\{ X\in M: M\models \varphi(X)\}\, .$$
Note that $$j(A)= \bigcup \{ j(A\cap V_\alpha): \alpha \in \OR\}$$
as $j(A\cap V_\alpha)=\{ X \in M: M\models \varphi(X)\}\cap V^M_{j(\alpha)}$.
Also note that if $A$ is a class of structures of the same type $\tau \in V_\kappa$,  then by elementarity $j(A)$ is also a subclass of $M$ of structures of  type $\tau$.

\subsection{$\Gamma_n$-strong cardinals}
\begin{definition} 
\label{defSigmaStrong}
For $n\geq 1$,  a cardinal $\kappa$  is \emph{$\lambda$-$\Gamma_n$-strong}  if for every $\Gamma_n$-definable (without parameters)  class $A$ there is   an elementary embedding $j:V\to M$, with $M$ transitive, $\crit(j)=\kappa$, $V_\lambda \subseteq M$,  and $A\cap V_\lambda \subseteq j(A)$.

$\kappa$ is \emph{$\Gamma_n$-strong} if it is $\lambda$-$\Gamma_n$-strong for every ordinal $\lambda$.

\end{definition}

Note that in the definition  above  $A \cap V_\lambda$ is only required to be contained in $j(A) \cap V_\lambda$ and not equal to it.  The reason is that in the $\Sigma_2$ case, if $A$ is the class of non-strong cardinals (which is $\Sigma_2$) and $\kappa$ is the least strong cardinal, then $\kappa \notin A$, but $\kappa \in j(A)$. See however the equivalence given in Proposition \ref{prop3.4}.

As with the case of strong cardinals, standard arguments show (cf. \cite{Kan:THI} 26.7(b)) that $\kappa$  is \emph{$\lambda$-$\Gamma_n$-strong} if and only if for every $\Gamma_n$-definable (without parameters)  class $A$ there is   an elementary embedding $j:V\to M$, with $M$ transitive, $\crit(j)=\kappa$, $V_\lambda \subseteq M$,  $j(\kappa)>\lambda$, and $A\cap V_\lambda \subseteq j(A)$.

\begin{proposition}
\label{prop12}
Every strong cardinal  is $\Sigma_2$-strong.
\end{proposition}

\begin{proof}
 Let $\kappa$ be a strong cardinal and let $A$ be a class that is $\Sigma_2$-definable (even allowing for parameters in $V_\kappa$). Let $\lambda \in C^{(2)}$ be greater than $\kappa$. Let $j:V\to M$ be elementary, with $M$ transitive, $\crit(j)=\kappa$, and $V_\lambda \subseteq M$. Let $\varphi$ be a $\Sigma_2$ formula defining $A$. If $a\in A\cap V_\lambda$, then $V_\lambda \models \varphi(a)$. Hence, since $V_\lambda \preceq_{\Sigma_1}M$, $M\models \varphi(a)$, and so $a\in j(A)=\{ x: M\models \varphi(x)\}$. 
\end{proof}

\begin{prop}
\label{prop3.3}
If $\lambda\in C^{(n+1)}$, then a cardinal $\kappa$ is $\lambda$-$\Pi_n$-strong if and only if it is $\lambda$-$\Sigma_{n+1}$-strong.
\end{prop}

\begin{proof}
Assume $\kappa$ is $\lambda$-$\Pi_n$-strong, with $\lambda\in C^{(n+1)}$, and let $A$ be a $\Sigma_{n+1}$-definable class. Let $\varphi (x)\equiv \exists y \psi (x,y)$ be a $\Sigma_{n+1}$ formula, with $\psi(x,y)$ being $\Pi_n$,  that defines $A$. Now define  $B$ as the class of all structures of the form $\langle V_\alpha , \in , a\rangle$, where $\alpha \in C^{(n)}$, $a\in V_\alpha$,  and $V_\alpha \models \varphi(a)$. Then $B$ is $\Pi_n$-definable. By our assumption, let $j:V\to M$ be an elementary embedding, with $M$ transitive, $\crit (j)=\kappa$, $V_\lambda \subseteq M$, and $B\cap V_\lambda \subseteq j(B)$. We just need to show that $A\cap V_\lambda \subseteq j(A)$. So, suppose $a\in A\cap V_\lambda$. Since $\lambda\in C^{(n+1)}$, we have that $V_\lambda \models \varphi(a)$. Let $b\in V_\lambda$ be a witness, so that $V_\lambda \models \psi(a,b)$. For some $\alpha <\lambda$ in $C^{(n)}$ we have that $a, b\in V_\alpha$. Hence, $V_\alpha \models \varphi(a)$. So $\langle V_\alpha , \in ,a\rangle \in B\cap V_\lambda$, and therefore  $\langle V_\alpha ,\in, a\rangle \in j(B)$. Thus, $M\models ``\alpha \in C^{(n)}, \, a\in V_\alpha, \mbox{ and }V_\alpha \models \varphi(a)"$. Hence, $M\models \varphi (a)$, i.e., $a\in j(A)$.
\end{proof}

\begin{corollary}
\label{coro3.4}
A cardinal $\kappa$ is $\Pi_n$-strong if and only if it is $\Sigma_{n+1}$-strong.
\end{corollary}

\begin{prop}
\label{prop3.3.1}
Suppose that $n\geq 2$ and $\lambda\in C^{(n)}$. Then the following are equivalent for a cardinal $\kappa <\lambda$:
\begin{enumerate}
\item $\kappa$ is $\lambda$-$\Sigma_{n}$-strong. 
\item There is an elementary embedding $j:V\to M$, with $M$ transitive, $\crit(j)=\kappa$, $V_\lambda \subseteq M$,  and  $M\models ``\lambda \in C^{(n-1)}"$.
\end{enumerate}
\end{prop}

\begin{proof}
(1)$\Rightarrow$(2): Suppose $\kappa$ is $\lambda$-$\Sigma_n$-strong. Let $A=C^{(n-1)}$. Since $A$ is $\Pi_{n-1}$-definable, hence also $\Sigma_n$-definable, by (1) there is an elementary embedding $j:V\to M$ with $M$ transitive, $\crit (j)=\kappa$, $V_\lambda \subseteq M$, and $A\cap V_\lambda \subseteq j(A)$. Since $\lambda \in C^{(n)}$, $C^{(n-1)}\cap \lambda$ is a club subset of $\lambda$. For every $\alpha <\lambda$ in $C^{(n-1)}$, $\alpha \in j(A)$, hence $M\models ``\alpha \in C^{(n-1)}"$ and so $M\models ``\lambda\mbox{ is a limit point of }C^{(n-1)}"$, which yields $M\models ``\lambda \in C^{(n-1)}"$.

(2)$\Rightarrow$(1): Let $A$ be a class definable by a $\Sigma_n$ formula $\varphi$,  and let $j:V\to M$ be an elementary embedding with $M$ transitive, $\crit(j)=\kappa$, $V_\lambda \subseteq M$,  and  $M\models ``\lambda \in C^{(n-1)}"$. Let $a\in A\cap V_\lambda$. Since $\lambda \in C^{(n)}$, $V_\lambda \models \varphi(a)$. And since $V_\lambda \subseteq M$ and $M\models ``\lambda \in C^{(n-1)}"$, $M\models \varphi(a)$, i.e., $a\in j(A)$. 
\end{proof}

The last proposition suggests the following definition and the ensuing characterization of $\Sigma_n$-strong cardinals in terms of extenders.

\begin{definition}
Given cardinals $\kappa <\lambda$, a \emph{$\Sigma_n$-strong $(\kappa,\lambda)$-extender} is a $(\kappa ,|V_{\lambda}|^+)$-extender $\Ee$ (see Definition \ref{defext})  such that $\overline{M}_\Ee\models ``\lambda \in C^{(n-1)}"$, where $\overline{M}_\Ee$ is the transitive collapse of the direct limit ultrapower $M_\Ee$  of $V$  by $\Ee$.
\end{definition}

\begin{prop}
\label{prop3.7}
If $n\geq 2$ and $\lambda\in C^{(n)}$, then a cardinal $\kappa <\lambda$ is $\lambda$-$\Sigma_n$-strong if and only if there exists a $\Sigma_n$-strong $(\kappa, \lambda)$-extender.
\end{prop}

\begin{proof}
If $\Ee$ is a $\Sigma_n$-strong $(\kappa ,\lambda)$-extender, then the extender embedding $j_\Ee:V\to \overline{M}_\Ee$,  witnesses that $\kappa$ is $\lambda$-$\Sigma_n$-strong (by Proposition \ref{prop3.3.1}).

Conversely, suppose  $j:V\to M$ is an elementary embedding, with $M$ transitive, $\crit(j)=\kappa$, $V_\lambda \subseteq M$,  and  $M\models ``\lambda \in C^{(n-1)}"$. Note that since $\lambda \in C^{(1)}$, $|V_\lambda|=\lambda$. Let $\Ee$ be the   $(\kappa , \lambda^+)$-extender derived from $j$. Namely,   for every $a\in [\lambda^+]^{<\omega}$ let $E_a$ be defined by:
$$X\in E_a \mbox{ if and only if } X\subseteq [\kappa]^{|a|} \mbox{ and }a \in j(X).$$
One can easily check that $\Ee$ satisfies conditions $(1)-(6)$ of Definition \ref{defext} (see \cite{Kan:THI} 26.7). So we only need to check that $\overline{M}_\Ee\models ``\lambda \in C^{(n-1)}"$. 

Let $j_\Ee:V\to M_\Ee$ and $k_\Ee:\overline{M}_\Ee\to M$ be the standard maps given by:
$j_\Ee(x)=[a, [c^a_x]]$ (any $a$), where $c^a_x :[\kappa]^{|a|} \to \{ x\}$; and $k_\Ee(\pi([a,[f]]))=j(f)(a)$, for $f:[\kappa]^{|a|}\to V$, where $\pi:M_\Ee\to \overline{M}_\Ee$ is the transitive collapse isomorphism. The maps $j_\Ee$ and $k_\Ee$ are elementary and $j=k_\Ee \circ j_\Ee$. Moreover, $k_\Ee\restriction V_\lambda$ is the identity. 

Since $M\models ``\lambda \in C^{(n-1)}"$, for each $\mu <\lambda$ in $C^{(n-1)}$, we have that $M\models ``\mu \in C^{(n-1)}"$. So, since $k_\Ee$ is elementary and is the identity on $V_\lambda$, we have that $M_\Ee\models ``\mu \in C^{(n-1)}"$. Hence, $M_\Ee\models ``\lambda \mbox{ is a limit point of }C^{(n-1)}"$, which yields $M_\Ee\models ``\lambda \in C^{(n-1)}"$.
\end{proof}

Similar characterizations may also be given for $\Pi_n$-strong cardinals. Notice that  (3) of the following proposition characterizes $\Pi_n$-strong cardinals as witnessing ``$\OR$ is Woodin"   restricted to $\Pi_n$-definable classes (see Definition \ref{defORWoodin}).

\begin{prop}
\label{prop3.4}
Suppose that $n\geq 1$ and $\lambda$ is a limit point of $C^{(n)}$. Then the following are equivalent for a cardinal $\kappa <\lambda$:
\begin{enumerate}
\item $\kappa$ is $\lambda$-$\Pi_{n}$-strong. 
\item There is an elementary embedding $j:V\to M$, with $M$ transitive, $\crit(j)=\kappa$, $V_\lambda \subseteq M$,  and  $M\models ``\lambda \in C^{(n)}"$.
\item For every $\Pi_n$-definable class $A$ there is an elementary embedding $j:V\to M$ with $M$ transitive, $\crit (j)=\kappa$, $V_\lambda \subseteq M$, and $A\cap V_\lambda = j(A)\cap V_\lambda$. 
\end{enumerate}
\end{prop}

\begin{proof}
(1)$\Rightarrow$(2): Suppose $\kappa$ is $\lambda$-$\Pi_n$-strong. Let $A=C^{(n)}$. Since $A$ is $\Pi_n$-definable, by (1)  there is an elementary embedding $j:V\to M$ with $M$ transitive, $\crit (j)=\kappa$, $V_\lambda \subseteq M$, and $A\cap V_\lambda \subseteq j(A)$. Thus, for every $\alpha <\lambda$ in $A$, $\alpha \in j(A)$, hence $M\models ``\alpha \in C^{(n)}"$ and so $M\models ``\lambda\mbox{ is a limit point of }C^{(n)}"$, which yields $M\models ``\lambda \in C^{(n)}"$.

(2)$\Rightarrow$(3): Let $A$ be a class definable by a $\Pi_n$ formula $\varphi (x)$,  and let $j:V\to M$ be an elementary embedding with $M$ transitive, $\crit(j)=\kappa$, $V_\lambda \subseteq M$,  and  $M\models ``\lambda \in C^{(n)}"$. Let $a\in A\cap V_\lambda$. Since $\lambda \in C^{(n)}$, $V_\lambda \models \varphi(a)$. And since $V_\lambda \subseteq M$ and $M\models ``\lambda \in C^{(n)}"$, $M\models \varphi(a)$, i.e., $a\in j(A)$. Conversely, suppose $a\in j(A)\cap V_\lambda$, i.e., $M\models ``\varphi(a)"$. Since $M\models ``\lambda \in C^{(n)}"$, $V_\lambda \models \varphi(a)$. And since $\lambda \in C^{(n)}$, $a\in A$. 

(3)$\Rightarrow$(1) is immediate.
\end{proof}

\begin{corollary}
\label{coroPinstrong}
Suppose that $n\geq 1$ and $\lambda$ is a limit point of $C^{(n)}$. A cardinal $\kappa$ is $\lambda$-$\Pi_n$-strong if and only if for every $\Pi_n$-definable class $A$ there is an elementary embedding $j:V\to M$ with $M$ transitive, $\crit (j)=\kappa$, $V_\lambda \subseteq M$, and $A\cap V_\lambda = j(A)\cap V_\lambda$.
\end{corollary}

Similarly as before we may also characterize $\Pi_n$-strong cardinals in terms of extenders. Namely,

\begin{definition}
Suppose that $n\geq 1$ and $\lambda$ is a limit point of $C^{(n)}$. Given a cardinal $\kappa <\lambda$, a \emph{$\Pi_n$-strong $(\kappa,\lambda)$-extender} is a $(\kappa ,|V_\lambda|^+)$-extender $\Ee$ (see Definition \ref{defext})  such that $\overline{M}_\Ee\models ``\lambda \in C^{(n)}"$, where $\overline{M}_\Ee$ is the transitive collapse of the direct limit ultrapower $M_\Ee$ of $V$ by $\Ee$.
\end{definition}

\begin{prop}
\label{prop3.10}
Suppose that $n\geq 1$ and $\lambda$ is a limit point of $C^{(n)}$. Then a cardinal $\kappa$ is $\lambda$-$\Pi_n$-strong if and only if there exists a $\Pi_n$-strong $(\kappa, \lambda)$-extender.
\end{prop}

It easily follows from the last proposition that being a $\Pi_n$-strong cardinal is a $\Pi_{n+1}$ property. Moreover, if $\kappa$ is $\Pi_n$-strong, then $\kappa \in C^{(n+1)}$. Hence, if $\kappa$ is $\Pi_{n+1}$-strong, then there are many $\Pi_n$-strong cardinals below $\kappa$, which shows that the $\Pi_n$-strong cardinals, $n>0$, form a hierarchy of strictly increasing strength. 

\medskip

Similarly as in \ref{prop1} we can prove the following.

\begin{prop}
\label{prop6}
If $\kappa$ is a $\Sigma_n$-strong cardinal, where $n\geq 2$, then $\Sigma_n(V_\kappa){\rm{-}}\PRP$ holds.
\end{prop}

\begin{proof}
Let $n\geq 2$. Let $\kappa$ be $\Sigma_n$-strong and let $\Ce$ be a  definable, by a $\Sigma_n$ formula  with parameters in $V_\kappa$,  proper class of structures in a fixed countable relational language. 
Let  $S:=\Ce\cap V_\kappa$. 

Given any $\Ae\in \Ce$, let $\lambda\geq \kappa$   with  $\Ae\in V_\lambda$.

Let $j:V\to M$ be an elementary embedding, with $crit(j)=\kappa$, $V_\lambda \subseteq M$,  $j(\kappa)>\lambda$, and $\Ce\cap V_\lambda \subseteq j(\Ce)$.

By elementarity, the restriction of $j$ to $\prod S$ yields a homomorphism
$$h:\prod S\to \prod (j(\Ce) \cap V^M_{j(\kappa)}).$$
Since  $\Ae\in \Ce \cap V_\lambda$, we have that  $\Ae \in j(\Ce)$.  Moreover $\Ae \in V_\lambda \subseteq V^M_{j(\kappa)}$. Thus, letting
$$g: \prod (j(\Ce) \cap V^M_{j(\kappa)})\to \Ae$$
be the projection map, we have that 
$$g\circ h:\prod S\to \Ae$$
is a homomorphism, as wanted.
\end{proof}

\subsection{The main theorem for $\Gamma_n$-strong cardinals}
Using similar arguments as in theorem \ref{thm1} we can now prove the main theorem of this section.

\begin{theorem}
\label{thm2}
The following are equivalent for $n\geq 2$:
\begin{enumerate}
\item There exists a  $\Sigma_n$-strong cardinal.
\item There exists a $\Pi_{n-1}$-strong cardinal.
\item $\Sigma_n$-$\PRP$
\item $\Pi_{n-1}$-$\PRP$
\item $\Sigma_n$-$\SWVP$
\item $\Pi_{n-1}$-$\SWVP$
\item $\Sigma_n$-$\WVP$
\item $\Pi_{n-1}$-$\WVP$
\end{enumerate}
\end{theorem}

\begin{proof}
(1)$\Rightarrow$(3) is given by proposition \ref{prop6}. (1)$\Leftrightarrow$(2) is given by corollary \ref{coro3.4}, (3)$\Rightarrow$(4),  (5)$\Rightarrow$(6), and (7)$\Rightarrow$(8) are immediate. The equivalences (3)$\Leftrightarrow$(5) and  (4)$\Leftrightarrow$(6) are given by corollary \ref{coro2}. The equivalence of (5) and (7), and also of (6) and (8), is given by theorem \ref{equivWandSW}. So, we only need to prove (4)$\Rightarrow$(2).

The proof is  analogous to the proof of Theorem  \ref{thm1}. So, we shall only indicate the relevant differences. Theorem \ref{thm1} proves the case $n=2$ (see proposition \ref{prop12}). 
Thus, we shall assume in the sequel that   $n> 2$.

\smallskip

Let $\mathcal{A}$ be the class of all structures  
$$\Ae_\alpha:=\langle V_{\alpha +1}, \in ,  \alpha , C^{(n-1)}\cap\alpha ,  \{ R^{\alpha}_\varphi \}_{\varphi \in \Pi_1} \rangle$$  where the constant $\alpha$  is the $\alpha$-th element of $C^{(n-1)}$,   and $\{ R^{\alpha}_\varphi\}_{\varphi \in \Pi_1}$ is the $\Pi_1$ relational diagram for $V_{\alpha +1}$, i.e., if $\varphi(x_1,\ldots ,x_n)$ is a $\Pi_1$ formula in the language of $\langle V_{\alpha +1},\in,\alpha, C^{(n-1)}\cap \alpha \rangle$, then $$R^\alpha_\varphi =\{ \langle x_1,\ldots ,x_n\rangle: \langle V_{\alpha +1},\in,\alpha, C^{(n-1)}\cap \alpha \rangle\models ``\varphi(x_1,\ldots ,x_n)"\}\, .$$

Then $\mathcal{A}$ is $\Pi_{n-1}$-definable without parameters. For $X\in \mathcal{A}$ if and only if $X=\langle X_0,X_1,X_2,X_3, X_4 \rangle$, where

 \begin{enumerate}
 \item[(1)] $X_2$  belongs to $C^{(n-1)}$ 
 \item[(2)] $X_0=V_{X_2 +1}$
 \item[(3)] $X_1 =\in \restriction X_0$
 \item[(4)] $X_0$ satisfies that $X_3=C^{(n-1)}\cap X_2$
 \item[(5)] $X_4$ is the $\Pi_1$ relational diagram of $\langle X_0,X_1,X_2,X_3 \rangle$
 \item[(6)] $\langle X_0, X_1, X_2, X_3 \rangle\models ``X_2$ is the $X_2$-th element of $C^{(n-1)}$".
\end{enumerate}

\medskip

Note that  the class $C$ of ordinals $\alpha$ such that $\Ae_\alpha \in \Ae$ is a closed and unbounded proper class.
By $\Pi_{n-1}$-$\PRP$  there exists a subset $S$ of $C$ such that for every $\beta\in C$ there is a homomorphism   $j_\beta : \prod_{\alpha \in S} \Ae_\alpha\to \Ae_\beta$. By enlarging $S$, if necessary, we may assume that $\kappa:= \rm{sup}(S)\in S$.

Now fix some $\beta\in C$  greater than $\kappa$,  of uncountable cofinality, and assume, towards a contradiction,  that no cardinal $\leq \kappa$ is $\beta$-$\Pi_{n-1}$-strong. Let $j=j_\beta$.

\medskip

From this point, the proof  proceeds  as in \ref{thm1}.  Namely, we define
$k:V_{\kappa +1}\to V_{\beta +1}$
by
$$k(X)=j(\langle X\cap V_{\alpha}\rangle_{\alpha \in S})$$
and note that $k(\kappa) =\beta$.

For each $a\in [\beta]^{<\omega}$,   define $E_a$ by
$$X\in E_a \quad \mbox{ iff }\quad X\subseteq [\kappa]^{|a|}  \mbox{ and } a\in k(X)\, .$$
As in \ref{thm1}, $E_a$ is an $\omega_1$-complete  ultrafilter over $[\kappa]^{|a|}$. Moreover, $\Ee :=\{ E_a: a\in [\beta]^{<\omega}\}$ is normal and coherent.

For each $a\in [\beta]^{<\omega}$, the ultrapower ${\rm{Ult}}(V,E_a)$  is well-founded, by $\omega_1$-completeness. So, let 
$$j_a:V\to M_a \cong {\rm{Ult}}(V, E_a)$$
with $M_a$ transitive, be the corresponding ultrapower embedding, and let    $M_\Ee$ be  the  direct limit of $$\langle \langle M_a:a\in [\beta]^{<\omega}\rangle, \langle i_{ab}:a\subseteq b\rangle\rangle$$ where the $i_{ab}:M_a\to M_b$ are the usual commuting  maps. The corresponding limit embedding $j_\Ee:V\to M_\Ee$ is elementary.
As in \ref{thm1}, 
$M_\Ee$ is closed under $\omega$-sequences, hence it is well-founded.  Moreover, letting $\pi:M_\Ee \to N$ be the transitive collapsing isomorphism, and  $j_N:V\to N$  the corresponding elementary embedding, i.e., $j_N=\pi \circ j_\Ee$, we have that $V_\beta \subseteq N$ and 
$j_N(\kappa)\geq \beta$.
Since $\beta >\kappa$, this implies that the critical point of $j_N$ is less than or equal to $\kappa$. The only additional argument needed, with respect to the proof of \ref{thm1},  is the following:

\begin{claim}
$N\models ``\beta \in C^{(n-1)}"$.
\end{claim}

\begin{proof}
Since $\beta$ is a limit point of $C^{(n-1)}$, it suffices to show that   if $\gamma <\beta$ belongs to  $C^{(n-1)}$, then $N\models ``\gamma \in C^{(n-1)}"$. So, fix some $\gamma <\beta$ in $C^{(n-1)}$.

Let $f :[\kappa]^1\to \kappa$ be such that $f (\{ x\})=x$. It is well known that $k_{\{ \gamma \}}( [f]_{E_{\{\gamma\}}})=\gamma$, where $k_{\{ \gamma\}}:M_{\{\gamma\}}\to N$ is the standard map given by $k_{\{ \gamma\}}([f]_{E_{\{\gamma\}}})=\pi([\{ \gamma\}, [f]_{E_{\{\gamma\}}}])$ (see \cite{Kan:THI} 26.2 (a)). 

Now notice that the set $X:=\{ \{ x\}\in [\kappa]^1: x\in C^{(n-1)}\}\in E_{\{\gamma\}}$, because $\{\gamma \}\in k(X)=\{\{ x\}\in [\beta]^1: x\in C^{(n-1)}\}$. Hence, $M_{\{ \gamma\}}\models ``[f]\in C^{(n-1)}"$, and therefore $M_\Ee \models ``[\{ \gamma\}, [f]]\in C^{(n-1)}"$, which yields $N\models ``\gamma \in C^{(n-1)}"$, as wanted. 
\end{proof}

Thus,   by proposition \ref{prop3.4}, $j_N$ witnesses that the critical point of $j_N$  is less than or equal to $\kappa$ and is $\beta$-$\Pi_{n-1}$-strong,  in contradiction to  our choice of $\beta$. 
\end{proof}

In a similar way we may obtain the following parameterized version of theorem \ref{thm2}. For  the proof of (4) implies (2),  we need to  consider the  class of structures 
$$\Ae_\alpha:=\langle V_{\alpha +1}, \in ,  \alpha , C^{(n-1)}\cap\alpha ,  \{ R^{\alpha}_\varphi \}_{\varphi \in \Pi_1} ,  \langle \delta \rangle_{\delta <\gamma} \rangle$$
where $$\Ae_\alpha:=\langle V_{\alpha +1}, \in ,  \alpha , C^{(n-1)}\cap\alpha ,  \{ R^{\alpha}_\varphi \}_{\varphi \in \Pi_1}\rangle$$ is as in the proof of Theorem \ref{thm2}, and we have a constant $\delta$ for every $\delta <\gamma$.

\begin{theorem}
\label{thm3}
The following are equivalent for $n\geq 2$:
\begin{enumerate}
\item There exist a  proper class of $\Sigma_n$-strong cardinals.
\item There exist a proper class of $\Pi_{n-1}$-strong cardinal.
\item $\mathbf{\Sigma_n}$-$\PRP$
\item $\mathbf{\Pi_{n-1}}$-$\PRP$
\item $\mathbf{\Sigma_n}$-$\SWVP$
\item $\mathbf{\Pi_{n-1}}$-$\SWVP$
\item $\mathbf{\Sigma_n}$-$\WVP$
\item $\mathbf{\Pi_{n-1}}$-$\WVP$
\end{enumerate}
\end{theorem}

\medskip

Recall that a cardinal $\kappa$ is \emph{Woodin} if for every  $A\subseteq V_\kappa$ there is $\alpha <\kappa$ such that $\alpha$ is $<\!\kappa$-$A$-strong, i.e., for every $\gamma <\kappa$ there is an elementary embedding $j:V\to M$ with $\crit (j)=\alpha$, $\gamma <j(\alpha)$, $V_\gamma \subseteq M$, and  $A\cap V_\gamma =j(A)\cap V_\gamma$. (See \cite{Kan:THI} 26.14.)

\begin{definition}
\label{defORWoodin}
\emph{$\OR$ is Woodin} if for every definable (with set parameters) class $A$ there exists some $\alpha$ which is $A$-strong, i.e., for every $\gamma$ there is an elementary embedding $j:V\to M$ with $\crit (j)=\alpha$, $\gamma <j(\alpha)$, $V_\gamma \subseteq M$, and  $A\cap V_\gamma =j(A)\cap V_\gamma$.
\end{definition}

The statement ``$\OR$ is Woodin" is first-order expressible as a schema, namely as ``There exists $\alpha$ which  is $A$-strong", for each definable $A$. Or equivalently, by corollary \ref{coroPinstrong},  as the schema ``There exists $\alpha$ which is $\Pi_n$-strong", for every $n$. Let us note that, by theorem \ref{thm3},  ``$\OR$ is Woodin" is also equivalent to the schema asserting ``There exist a proper class of $\alpha$ which are $\Pi_n$-strong", for every $n$.  Thus, theorem \ref{thm3} yields the following corollary, first proved by the second author in \cite{W:LCSWVP}, which gives the exact large-cardinal strength of $\WVP$ and $\SWVP$.

\begin{corollary}
\label{maincoro}
The following are equivalent:
\begin{enumerate}
\item $\OR$ is Woodin 
\item $\SWVP$
\item $\WVP$.
\end{enumerate}
\end{corollary}

\end{document}